\documentclass[12pt,ruled,noline]{article}
\usepackage[utf8]{inputenc}
\usepackage[a4paper]{geometry}
\usepackage{color}
\usepackage{array}
\usepackage{booktabs}
\usepackage{textcomp}
\usepackage{mathtools}
\usepackage{multirow}
\usepackage{algorithm2e}
\usepackage{amsmath}
\usepackage{amsthm}
\usepackage{amssymb}
\usepackage{graphicx}
\usepackage[unicode=true,
 bookmarks=false,
breaklinks=false,pdfborder={0 0 1},backref=false,colorlinks=true]
{hyperref}
\hypersetup{
 linkcolor=blue,citecolor=darkgreen,urlcolor=blue}
 
 \usepackage{float}
 
\definecolor{darkgreen}{RGB}{0,150,0}
\makeatletter

\providecommand{\tabularnewline}{\\}

\theoremstyle{plain}
\newtheorem{thm}{\protect\theoremname}
\theoremstyle{definition}
\newtheorem{problem}[thm]{\protect\problemname}
\theoremstyle{definition}
\newtheorem{example}[thm]{\protect\examplename}
\theoremstyle{definition}
\newtheorem{defn}[thm]{\protect\definitionname}
\theoremstyle{plain}
\newtheorem{fact}[thm]{\protect\factname}

\usepackage{amsfonts}\usepackage{enumerate}

\usepackage{graphicx}
\SetArgSty{textup}
\usepackage{array}
\usepackage{subcaption}
\usepackage{todonotes}
\usepackage{dcolumn}

\LinesNumbered

\SetKw{Set}{set}{}
\SetKwBlock{Repeat}{repeat}{}

\DeclareMathOperator*{\argmin}{argmin}

\makeatletter
\let\@fnsymbol\@arabic
\makeatother

\makeatother

\providecommand{\definitionname}{Definition}
\providecommand{\examplename}{Example}
\providecommand{\factname}{Fact}
\providecommand{\problemname}{Problem}
\providecommand{\theoremname}{Theorem}

\begin{document}
\title{\textbf{The superiorization method with restarted perturbations for
split minimization problems with an application to radiotherapy treatment
planning}}
\author{Francisco J.\ Aragón-Artacho\thanks{Department of Mathematics, University of Alicante, Alicante, \textsc{Spain}.
\protect \\
Email:~\protect\href{mailto:francisco.aragon@ua.es}{francisco.aragon@ua.es}} \and Yair Censor\thanks{Department of Mathematics, University of Haifa, Haifa, \textsc{Israel}.
\protect \\
Email: \protect\href{mailto:yair@math.haifa.ac.il}{yair@math.haifa.ac.il}} \and Aviv Gibali\thanks{Department of Mathematics, Braude College, Karmiel, \textsc{Israel}.
\protect \\
Email: \protect\href{mailto:avivg@braude.ac.il}{avivg@braude.ac.il}} \and David Torregrosa-Belén\thanks{Department of Mathematics, University of Alicante, Alicante, \textsc{Spain}.
\protect \\
Email: \protect\href{mailto:david.torregrosa@ua.es}{david.torregrosa@ua.es}}}
\date{July 12, 2022. Revised: October 10, 2022}
\maketitle
\begin{abstract}
In this paper we study the split minimization problem that consists
of two constrained minimization problems in two separate spaces that
are connected via a linear operator that maps one space into the other.
To handle the data of such a problem we develop a superiorization
approach that can reach a feasible point with reduced (not necessarily
minimal) objective function values. The superiorization methodology
is based on interlacing the iterative steps of two separate and independent
iterative processes by perturbing the iterates of one process according
to the steps dictated by the other process. We include in our developed
method two novel elements. The first one is the permission to restart
the perturbations in the superiorized algorithm which results in a
significant acceleration and increases the computational efficiency.
The second element is the ability to independently superiorize subvectors.
This caters to the needs of real-world applications, as demonstrated
here for a problem in intensity-modulated radiation therapy treatment
planning.
\end{abstract}

\section{Introduction}

In a fair number of applications the nature and size of the arising
constrained optimization problems make it computationally difficult,
or sometimes even impossible, to obtain exact solutions and alternative
ways of handling the data of the optimization problem should be considered.
A common approach is the regularization technique that replaces the
constrained optimization problem by an unconstrained optimization
problem wherein the objective function is a linear combination of
the original objective and a regularizing term that ``measures''
in some way the constraints violations.

This approach is used for constrained minimization problems appearing
in image processing, where the celebrated Fast Iterative Shrinkage-Thresholding
Algorithm (FISTA) method was pioneered by Beck and Teboulle \cite{FISTA}.
In situations when there are some constraints whose satisfaction is
imperative (``hard constraints'') the problem can be considered
as being composed of two goals: A major goal of satisfying the constraints
and a secondary, but desirable, goal of target (a.k.a. objective,
merit, cost) function reduction, not necessarily minimization.

In this setting, the \emph{superiorization methodology} (SM) has proven
capable of efficiently handling the data of very large constrained
optimization problems. The idea behind superiorization is to apply
a feasibility-seeking algorithm and introduce in each of its iterations
a certain change, referred to as a \emph{perturbation}, whose aim
is to reduce the value of the target function. When the feasibility-seeking
iterative algorithm is \emph{bounded perturbation resilient} (see
Definition~\ref{def:bpr} below), the superiorized version of the
feasibility-seeking algorithm will converge to a feasible solution
which is expected to have a reduced, not necessarily minimal, target
function value.

We elaborate on the SM in Section \ref{Sect:2} below. The reader
might find the online bibliography \cite{sup-bib-online}, which is
an updated snapshot of current work in this field, useful. In particular,
\cite{weak-2015,Humphris-2021,herman-2020}, to name but a few, and
the references therein, contain relevant introductions to the SM.

Superiorization has successfully found multiple applications, in some
cases outperforming other state-of-the-art algorithms, such as computed
tomography~\cite{supvssubg,WHare}, inverse treatment
planning in radiation therapy~\cite{sup:radthe}, bioluminiscence
tomography~\cite{sup:bio}  and linear optimization~\cite{LinSup}.
Many more are documented in~\cite{sup-bib-online}.

Regarding the question how SM-based algorithms compare with exact
constrained optimization algorithms, it should be noted that the SM
is not intended to solve exactly constrained optimization problems, in spite of that such comparisons have been reported. There is an
extensive comparison between superiorization methods and regularization
methods in \cite{WHare}, see also \cite{censor-petra-2020}. In
\cite{supvssubg} a comparison of SM with the projected subgradient
method shows, surprisingly, the advantages of SM within the computed
tomography application discussed there.

The paper \cite{LinSup} compares the performance of a linear superiorized
method (LinSup) and the Simplex algorithm built-in Matlab's \texttt{linprog}
solver, for linear optimization problems. The numerical experiments
there show that for large-scale problems, the use of LinSup can be
advantageous to the employment of the Simplex algorithm.

In this paper we propose a novel superiorized algorithm for dealing
with the data of the following \emph{split minimization problem:}
\begin{problem}
\label{prob:smp-1}\textbf{\emph{The split minimization problem}}\textbf{
(SMP)}. Given two nonempty, closed and convex subsets $C\subseteq\mathbb{R}^{n}$
and $Q\subseteq\mathbb{R}^{m}$ of two Euclidean spaces of dimensions
$n$ and $m$, respectively, an $m\times n$ real matrix $A$, and
convex functions $f:\mathbb{R}^{n}\to\mathbb{R}$ and $g:\mathbb{R}^{m}\to\mathbb{R},$
find
\begin{equation}
\begin{aligned}x^{*}\in C\;\text{such that}\;x^{*}\in\argmin\{f(x)\mid x\in C\}\end{aligned}
,\;\text{ and such that}\;\label{e:svip1}
\end{equation}
\begin{equation}
\begin{aligned}y^{*}:=Ax^{*}\in Q\;\text{and}\;y^{*}\in\argmin\{g(y)\mid y\in Q\}.\end{aligned}
\label{e:svip2}
\end{equation}
\end{problem}

It is important to observe that the two objective functions $f$ and $g$ in~\eqref{e:svip1} and~\eqref{e:svip2} may conflict with each other, and thus the existence of a solution to Problem~\ref{prob:smp-1} is not guaranteed even if $Ax\in Q$ for all $x\in C$. This is a new genre of problems which
are not considered as multi-objective but rather split between two spaces. Problem~\ref{prob:smp-1} is a particular instance of the \emph{split
variational inequality problem} (SVIP), which employs, instead of
the minimization problems in~\eqref{e:svip1} and~\eqref{e:svip2},
variational inequalities. The SVIP, see~\cite{svip}, entails
finding a solution of one \textit{variational inequality problem}
(VIP), the image of which under a given bounded linear transformation
is a solution of another VIP. Algorithms
for solving the SVIP require computing the projections onto the corresponding
constraint sets at every step, see~\cite{svip}. In the case when
$C$ and $Q$ are each given by an intersection of
nonempty, closed and convex sets, auxiliary algorithms, such as
Dykstra's algorithm~\cite{dijkstra} (see also \cite[Subsection 30.2]{bauschke2017}),
the Halpern--Lions--Wittmann--Bauschke (HLWB) algorithm~\cite{HLWB}
(see also \cite[Subsection 30.1]{bauschke2017}) or the averaged
alternating modified reflections method~\cite{aamr} are needed
for computing/approximating these projections, which will considerably
increase the running time and the numerical errors of the algorithms.
  In this work, we do not aim to find an exact solution of the SMP, but rather obtain a feasible solution with reduced values of the objective functions $f$ and $g$. This allows us to drop the usual assumption on the existence of a solution to the minimization problem, and instead we will only require that the set of solutions to the associated feasibility problem (see, Problem~\ref{prob:mssfp}) is nonempty.

Two novel elements are included here in our superiorized algorithm
for the SMP. The first is a permission to restart the perturbations
in the superiorized algorithm which increases the computational efficiency.
The second is the ability to superiorize independently over subvectors.
This caters to real-world situations, as we demonstrate here for a
problem in \emph{intensity-modulated radiation therapy} (IMRT) treatment
planning.

The remainder of the paper is structured as follows. In Section \ref{Sect:2}
we briefly present the superiorization methodology and in Section
\ref{sec:The-superiorization-method}, we introduce a new technique
for setting up the step-sizes in the perturbations, which results
in a new version of the general structure of the superiorized algorithm.
This new structure increases the efficiency of the superiorized algorithm
by allowing restarts of step-sizes. Our new algorithm for dealing
with the data of the split minimization problem is developed in Section~\ref{Sect:3}.
Finally, in Section~\ref{Sect:4} we present some numerical experiments
on three demonstrative examples of the performance of the algorithms.
The last example is a nontrivial realistic problem arising in IMRT.

\section{The superiorization methodology}

\label{Sect:2}

In this section we present a brief introduction to the superiorization
methodology (SM)\footnote{A word about the history: The superiorization method was born when
the terms and notions “superiorization” and “perturbation resilience”,
in the present context, first appeared in the 2009 paper of Davidi,
Herman and Censor \cite{bip_bpr} which followed its 2007 forerunner
by Butnariu et al. \cite{sm2007}. The ideas have some of their roots
in the 2006 and 2008 papers of Butnariu et al. \cite{brz06} and \cite{brz08}.
All these culminated in Ran Davidi’s 2010 PhD dissertation \cite{davidi-phd}
and the many papers since then cited in \cite{sup-bib-online}, such
as, e.g., \cite{Humphris-2021}.}, which is a simplified version of the presentation in~\cite{SupCDH}.
The SM has been shown to be a useful tool for handling the data of
difficult constrained minimization problems of the form 
\begin{equation}
\min\{\phi(x)\mid x\in C\},\label{e:cm}
\end{equation}
where $\phi:\mathbb{R}^{n}\to\mathbb{R}$ is a target function and
$C\subseteq\mathbb{R}^{n}$ is a nonempty feasible set, generally
presented as an intersection of a finite family of constraint sets
$C:=\cap_{s=1}^{p}C_{s}$. When $\{C_{s}\}_{s=1}^{p}$ is a family
of nonempty, closed and convex sets in $\mathbb{R}^{n}$, there is
a wide range of projection methods (see, e.g.,~\cite{bauschke2017,cegielski})
that can be employed for solving the convex feasibility-seeking problem
\begin{equation}
\text{find}\;x^{*}\in C:=\cap_{s=1}^{p}C_{s}.\label{e:fp}
\end{equation}
The first building brick of the SM is an iterative feasibility-seeking
algorithm, often a projection method, which is referred to as the
\emph{basic algorithm}, capable of finding (asymptotically) a solution
to~\eqref{e:fp}. This algorithm employs an \emph{algorithmic operator}
$\mathbf{T}_{C}:\mathbb{R}^{n}\to\mathbb{R}^{n}$ in the following
iterative process.

\smallskip{}

\begin{algorithm}
\textbf{Initialization:}\label{alg:basic} Choose an arbitrary initialization
point $x^{0}\in\mathbb{R}^{n}$\; \textbf{Iterative Step:}~Given
the current iterate $x^{k}$, calculate the next iterate $x^{k+1}$
by 
\begin{equation}
x^{k+1}=\mathbf{T}_{C}(x^{k}).\label{e:basic}
\end{equation}
\caption{The Basic Algorithm}
\end{algorithm}

\smallskip{}

\begin{example}
A well-known feasibility-seeking algorithm for the set $C$ given
in~\eqref{e:fp} is the method of \emph{sequential alternating projections},
whose algorithmic operator is given by 
\begin{equation}
\mathbf{T}_{C}:=P_{C_{p}}P_{C_{p-1}}\cdots P_{C_{2}}P_{C_{1}},\label{e:ap}
\end{equation}
where $P_{C_{s}}$ denotes the \emph{orthogonal projection }onto the
set $C_{s}$, i.e., 
\begin{equation}
P_{C_{s}}(x):=\argmin\{\|x-c\|\mid c\in C_{s}\}.\label{e:projection}
\end{equation}
\end{example}

Many other iterative feasibility-seeking projection methods are available,
see, e.g., the excellent review paper of Bauschke and Borwein \cite{BauschkeBorwein1996jour}
and \cite{escalante-2011}. Such methods have general algorithmic
structures of block-iterative projection (BIP), see, e.g., \cite{AharoniCensor1989jour,AleynerReich2008jour}
or string-averaging projections (SAP), see, e.g., \cite{sap-pink-2001,nisenbaum,Nikazad-sap-2016}.

In the SM one constructs from the basic algorithm a ``superiorized
version of the basic algorithm'' which includes perturbations of
the iterates of the basic algorithm. This requires the basic algorithm
to be resilient to certain perturbations. The definition is given
next with respect to the feasibility-seeking operator \eqref{e:ap},
but is phrased in the literature with algorithmic operators of any
basic algorithm.
\begin{defn}
\textbf{\textit{Bounded perturbation resilience}}. \label{def:bpr}
Let $\left\{ C_{s}\right\} _{s=1}^{p}$ be a family
of closed and convex sets in $\mathbb{R}^{n}$ such that $C=\bigcap_{s=1}^{p}C_{s}$
is nonempty. An algorithmic operator $\mathbf{T}_{C}:\mathbb{R}^{n}\to\mathbb{R}^{n}$
for solving the feasibility-seeking problem associated with $C$ is
said to be \emph{bounded perturbation resilient} if the following
holds: for all $x^{0}\in\mathbb{R}^{n}$, if the sequence $\{x^{k}\}_{k=0}^{\infty}$
generated by Algorithm~\ref{alg:basic} converges to a solution of
the feasibility-seeking problem, then any sequence $\{y^{k}\}_{k=0}^{\infty}$
generated by 
\begin{equation}
y^{k+1}=\mathbf{T}_{C}\left(y^{k}+\eta_{k}v^{k}\right),\quad\text{for all }k\geq0,\label{e:perturb}
\end{equation}
for any $y^{0}\in\mathbb{R}^{n},$ where the vector sequence $\{v^{k}\}_{k=0}^{\infty}$
is bounded and the scalars $\{\eta_{k}\}_{k=0}^{\infty}$ are nonnegative
and summable, i.e., ${\textstyle \sum_{k=0}^{\infty}\eta_{k}<\infty,}$
also converges to a solution of the feasibility-seeking problem.
\end{defn}

The property of bounded perturbation resilience has been validated
for two major prototypical algorithmic operators that give rise to
the \emph{string averaging projections} method and the \emph{block
iterative projections} method~mentioned above, see~\cite{bip_bpr}
and~\cite{sm2007}, respectively. These schemes include many well-known
projection algorithms, such as the method of alternating projections
and Cimmino's algorithm. The convexity and closedness of the sets
$C_{s}$ is present in the results proving the bounded perturbation
resilience of the BIP algorithms and the SAP methods.

The importance of bounded perturbation resilience for the SM stems
from the fact that it allows to include perturbations in the iterative
steps of the basic algorithm without compromising its convergence
to a feasible solution, while steering the algorithm toward a feasible
point with a reduced (not necessarily minimal) value of the target
function.

The fundamental idea underlying the SM is to use the bounded perturbations
in~\eqref{e:perturb} in order to induce convergence to a feasible
point which is \emph{superior, }meaning that the value of the target
function $\phi$ is smaller or equal than that of a point obtained
by applying the basic algorithm alone without perturbations. To achieve
this aim, the bounded perturbations in~\eqref{e:perturb} should
imply that 
\begin{equation}
\phi(y^{k}+\eta_{k}v^{k})\leq\phi(y^{k}),\quad\text{for all }k\geq0.\label{e:supphi}
\end{equation}
To do so, the sequence $\{v^{k}\}_{k=0}^{\infty}$ is chosen according
to the next definition, which is closely related to the concept of
\emph{descent direction}.
\begin{defn}
Given a function $\phi:\mathbb{R}^{n}\to\mathbb{R}$ and a point $y\in\mathbb{R}^{n}$,
we say that a direction $v\in\mathbb{R}^{n}$ is \emph{nonascending
for} $\phi$ \emph{at} $y$ if $\|v\|\leq1$ and there is some $\delta>0$
such that 
\begin{equation}
\phi(y+\lambda v)\leq\phi(y),\quad\text{for all }\lambda\in{[0,\delta]}.\label{e:nonascd}
\end{equation}
\end{defn}

Obviously, the zero vector is a nonascending direction. However, it
would not provide any perturbation of the sequence in~\eqref{e:perturb}.
Denoting by $\frac{\partial\phi}{\partial x_{i}}(x)$ the partial derivatives,
the next result provides a formula for obtaining nonascending vectors
of convex functions.
\begin{thm}
\label{t:nonasc}\cite[Theorem~2]{SM2012}. Let $\phi:\mathbb{R}^{n}\to\mathbb{R}$
be a convex function and let $x\in\mathbb{R}^{n}$. Let $u=(u_{i})_{i=1}^{n}\in\mathbb{R}^{n}$
be defined, by
\begin{equation}
u_{i}:=\left\{ \begin{array}{ll}
\frac{\partial\phi}{\partial x_{i}}(x), & \text{if }\frac{\partial\phi}{\partial x_{i}}(x)\text{ exists},\\
0, & \text{otherwise},
\end{array}\right.\label{e:g}
\end{equation}
and define 
\begin{equation}
v:=\left\{ \begin{array}{ll}
0, & \text{if }\|u\|=0,\\
-u/\|u\|, & \text{otherwise}.
\end{array}\right.
\end{equation}
Then $v$ is a nonascending vector for $\phi$ at the point $x$.
\end{thm}

Next we present the pseudo-code of the iterative process governing
the superiorized version of the basic algorithm.

\smallskip{}

\begin{algorithm}
\textbf{Initialization:}~Choose an arbitrary initialization point
$y^{0}\in\mathbb{R}^{n}$, a summable nonnegative sequence $\{\eta_{\ell}\}_{\ell=0}^{\infty}$
and a positive integer $N$\; \Set{$k=0$ and $\ell=-1$}\; \Repeat{
\Set{$y^{k,0}=y^{k}$}\; \For{$j=0$ \KwTo $N-1$}{ \Set{$v^{k,j}$
to be a nonascending vector for $\phi$ at $y^{k,j}$}\;\Set{$\ell=\ell+1$}\;
\While{$\phi(y^{k,j}+\eta_{\ell}v^{k,j})>\phi(y^{k})$}{ \Set{$\ell=\ell+1$}\;
} \Set{$y^{k,j+1}=y^{k,j}+\eta_{\ell}v^{k,j}$}\; } \Set{$y^{k+1}=\mathbf{T}_{C}\left(y^{k,N}\right)$
and $k=k+1$}\; } \caption{The Superiorized Version of the Basic Algorithm}
\label{a:SV}
\end{algorithm}

\smallskip{}
 The choice of nonascending vectors guarantees that the \textbf{while}
loop in lines 8--10 is finite (see \cite{SM2012} for a complete
proof on the termination of the algorithm). When a bounded perturbation
resilient operator $\mathbf{T}_{C}$ is chosen as the basic algorithm,
Algorithm~\ref{a:SV} will converge to a solution of the feasibility-seeking
problem. Moreover, it is expected that the perturbations in line 11,
which reduce at each inner-loop step the value of the target function
$\phi$ (line~8), will drive the iterates of Algorithm~\ref{a:SV}
to an output which will be superior (from the point of view of its
$\phi$ value) to the output that would have been obtained by the
original unperturbed basic algorithm.

\subsection{The guarantee problem of the SM \label{rem:failure}}

Proving mathematically a guarantee of global function reduction of
the SM will probably require some additional assumptions on the feasible
set, the target function, the parameters involved, or even on the
initialization points. We present here a simple example where the
performance of the SM depends on the choice of the initialization
point.

\begin{figure}[h]
\centering{}\includegraphics[width=0.45\textwidth,height=0.315\textwidth]{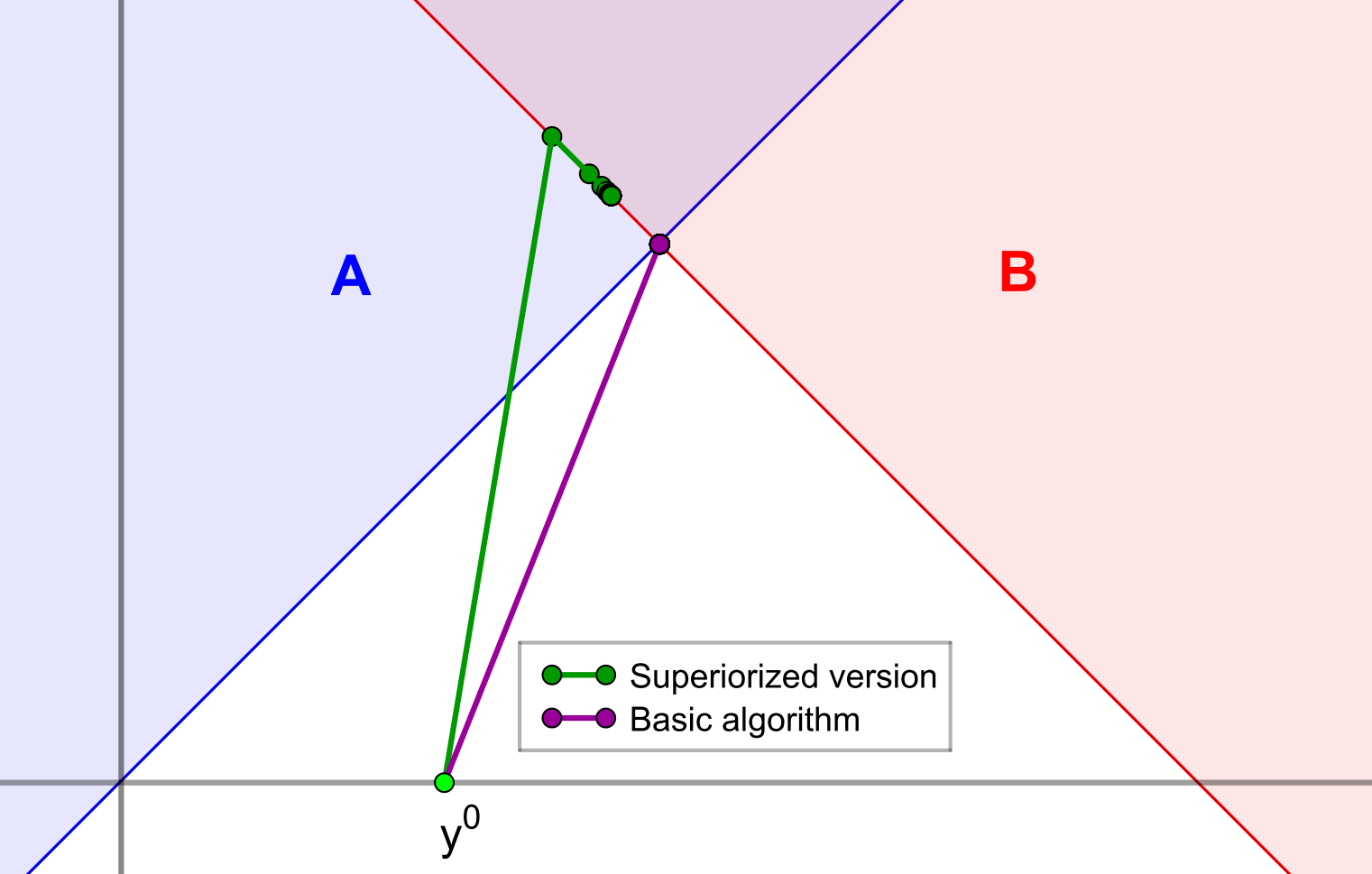}
\hspace{0.01\textwidth}\includegraphics[width=0.45\textwidth,height=0.315\textwidth]{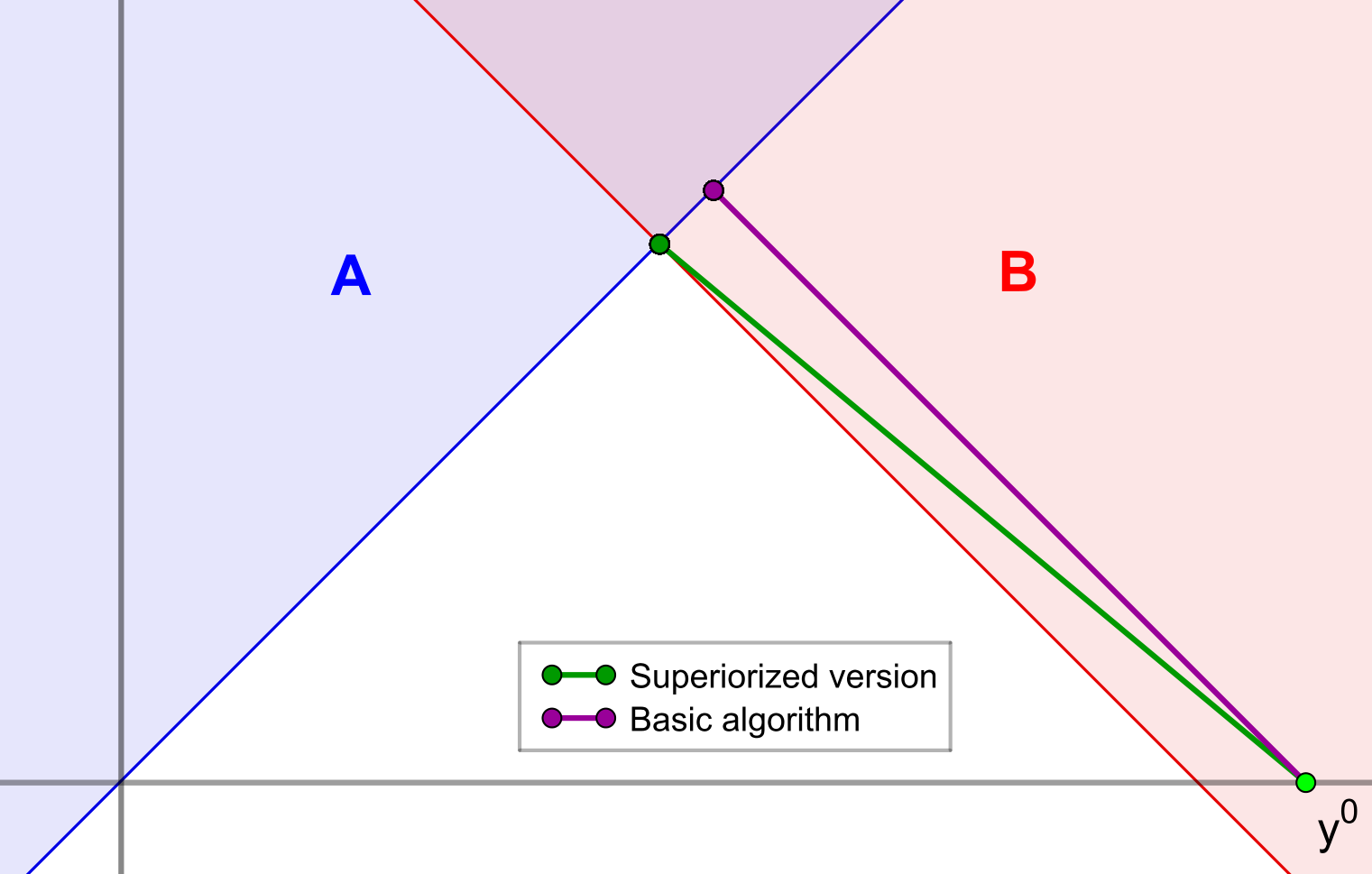}\caption{(Left) Taking $y^{0}=(3/10,0)^{T}$ as starting point, the alternating
projection method converges in one iteration to the minimum norm point
in the intersection of two half-spaces, while its superiorized version
remains far from the solution after 50 iterations. (Right) If we
take $y^{0}=(11/10,0)^{T},$ in one iteration, the superiorized version
converges to the minimum norm point while the alternating projection
method converges to a feasible point which is not the solution. \label{fig:failure}}
\end{figure}

Consider the problem of finding the minimum norm point in the intersection
of two half-spaces $A:=\left\{ x\in\mathbb{R}^{2}\mid x_{1}+x_{2}\geq1\right\} $
and $B:=\left\{ x\in\mathbb{R}^{2}\mid x_{1}-x_{2}\leq0\right\} $.
If we use the method of alternating projections as the basic algorithm
with starting point $y^{0}:=(3/10,0)^{T}$, one obtains $y^{1}=(1/2,1/2)^{T}$,
which is the solution to the problem. In Figure~\ref{fig:failure}\footnote{All figures in this paper, including this one, are created by the
Cinderella interactive geometry software \cite{cinderella}.} (left) we show 50 iterations of its superiorized version for the
same starting point with $\phi(y):=\|y\|^{2}$, step-sizes in the
sequence $\left\{ 2^{-\ell}\right\} _{\ell=1}^{\infty}$, $N=1$ and
taking as nonascending direction $v^{k}:=-\frac{\nabla\phi(y^{k})}{\|\nabla\phi(y^{k})\|}=-\frac{y^{k}}{\|y^{k}\|}$.

From this starting point the basic algorithm for feasibility-seeking
alone without perturbations (i.e., the alternating projection method)
converges in one iteration to the minimum norm point in the intersection
of the two half-spaces, while its superiorized version remains far
from the solution after 50 iterations. This happens because the
first perturbation applied to this $y^{0}$ results in a point on
the horizontal axis inside the set $A.$ If we were to choose another
$y^{0}$ on the $x$-axis, but far enough to the right inside the
set $B$, then after one perturbation the next point will be outside
both sets $A$ and $B$ on the positive $x$-axis and this would lead,
following a single iteration of feasibility-seeking, to the minimum
norm point whereas the feasibility-seeking-only from that $y^{0}$
onward would lead to a less ``superior'' feasible point, as shown
in Figure \ref{fig:failure} (right).

Observe that we only computed 50 iterations because after them the norm of the perturbations is smaller than $2^{-50}\approx8.9\cdot10^{-16}$.
Hence, the effect of the perturbations steering the basic algorithm
to a superiorized solution vanishes, having no real effect on it after
even less than 100 iterations. This phenomenon is inherent in the
SM and the purpose of the restarts, proposed in the next section,
is to improve this unwanted behavior.

In \cite[Section 3]{censor-levy-2019} we gave a precise
definition of the ``guarantee problem'' of the SM. We wrote there:
``The SM interlaces into a feasibility-seeking basic algorithm target
function reduction steps. These steps cause the target function to
reach lower values locally, prior to performing the next feasibility-seeking
iterations. A mathematical guarantee has not been found to date that
the overall process of the superiorized version of the basic algorithm
will not only retain its feasibility-seeking nature but also accumulate
and preserve globally the target function reductions. We call this
fundamental question of the SM ``the guarantee problem of the SM''\ which
is: under which conditions one can guarantee that a superiorized version
of a bounded perturbation resilient feasibility-seeking algorithm
converges to a feasible point that has target function value smaller
or equal to that of a point to which this algorithm would have converged
if no perturbations were applied -- everything else being equal.''\

Numerous works that are cited in \cite{sup-bib-online} show that
this global function reduction of the SM occurs in practice in many
real-world applications. In addition to a partial answer in \cite{censor-levy-2019}
with the aid of the ``concentration of measure'' principle there
is also the partial result of \cite[Theorem 4.1]{cz3-2015} about
strict Fejér monotonicity of sequences generated by an SM algorithm.

\section{First development: The superiorization method with perturbations
restarts}

\label{sec:The-superiorization-method}

We offer here a modification of the SM, applied to the superiorized
version of the basic algorithm in Algorithm~\ref{a:SV}, by setting
the perturbations step-sizes in a manner that allows restarts. Commonly,
the summable sequence $\{\eta_{\ell}\}_{\ell=0}^{\infty}$ employed
in Algorithm~\ref{a:SV} is being generated by taking a real number
$\alpha\in{]0,1[}$, referred to as \emph{kernel}, and setting $\eta_{\ell}:=\alpha^{\ell}$,
for $\ell\geq0$. This strategy works well in practice, as witnessed
by many works cited in \cite{sup-bib-online}, but it has though the
inconvenience that, as the iterations progress, the step-sizes in~\eqref{e:perturb}
decrease toward zero quite fast, yielding insignificant perturbations.

In some applications, various methods have been studied for controlling
the step-sizes, see, e.g.,~\cite{LinSup,SupStep1,SupStep2}, see also the software package SNARK14 \cite{SNARK} which is an updated version of \cite{SNARKpre}. We
propose here a new strategy which allows \emph{restarting} the sequence
of step-sizes to a previous value while maintaining the summability
of the series of step-sizes. The restarting of step-sizes is a useful
approach that allows to improve an algorithm's performance, see, e.g.
\cite{restartSGD,restFista}, where this technique is applied to
the stochastic gradient descent and FISTA, respectively.

Our proposed scheme for restarting the step-sizes is controlled by
a sequence of positive integers $\{W_{r}\}_{r=0}^{\infty}$ and we
call the indices $r=1,2,\ldots$ \textit{restart indices}. Bearing
some similarity to a backtracking scheme, the initial step-size at
the beginning of a new loop of perturbations is reduced after each
restart. Specifically, let $\alpha\in{]0,1[}$ be any fixed kernel,
assume that $r$ restarts have been already performed and let $W_{r}$
denote the number of consecutive step-sizes which must be taken before
allowing the next restart.

The algorithm begins with $r=0$ and takes $W_{0}$ decreasing step-sizes
in the sequence $\{1,\alpha,\alpha^{2},\ldots\}$. After these $W_{0}$
step-sizes, the algorithm performs a restart in the step-sizes by
setting $r=1$ and taking anew $W_{1}$ decreasing step-sizes in the
sequence $\{\alpha,\alpha^{2},\alpha^{3},\ldots\}$. Then the algorithm
performs another restart with $r=2$ and uses step-sizes in the sequence
$\{\alpha^{2},\alpha^{3},\alpha^{4},\ldots\}$, and so on.

This is accurately described in the pseudo-code of the \textit{Superiorized
Version of the Basic Algorithm with Restarts} presented below. Observe
that there are no restrictions on how the sequence $\{W_{r}\}_{r=0}^{\infty}$
is chosen. A simple possible choice is to take a positive constant
value $W_{r}:=W$, for all $r\geq0$. 

Note that, since the kernel $\alpha$ needs to lie
in the interval $]0,1[$, despite performing restarts, the perturbations
may happen to yield an insignificant decrease in the objective value,
specially if the norm of the nonascending vectors is close to zero.
This can be controlled by considering a positive real number $c$
and performing the restarts to the sequence $\{c\,\alpha^{\ell}\}_{\ell=0}^{\infty}$
rather than to $\{\alpha^{\ell}\}_{\ell=0}^{\infty}$. This parameter
is also included in the pseudo-code of Algorithm~\ref{a:RSV-1}.

\begin{algorithm}[ht!]
\textbf{Initialization:}~Choose an arbitrary initialization point
$y^{0}\in\mathbb{R}^{n}$, $\alpha\in{]0,1[}$, a
positive number $c$, a positive integer $N$ and a sequence of
positive integers $\{W_{r}\}_{r=0}^{\infty}$ \; \Set{$k=0$, $\ell=-1$,
$w=0$ and $r=0$}\; \Repeat{\Set{$y^{k,0}=y^{k}$}\; \For{$j=0$
\KwTo $N-1$}{ \Set{$v^{k,j}$ to be a nonascending vector for
$\phi$ at $y^{k,j}$}\;\Set{$\ell=\ell+1$}\; \While{$\phi(y^{k,j}+c\,\alpha^{\ell}v^{k,j})>\phi(y^{k})$}{
\Set{$\ell=\ell+1$}\; } \Set{$y^{k,j+1}=y^{k,j}+c\,\alpha^{\ell}v^{k,j}$}\;
} \Set{$w=w+1$}\; \If{$w=W_{r}$}{ \Set{$r=r+1$, $\ell=r$
and $w=0$}\; } \Set{$y^{k+1}=\mathbf{T}_{C}\left(y^{k,N}\right)$
and $k=k+1$}\; } \caption{The Superiorized Version of the Basic Algorithm with Restarts}
\label{a:RSV-1}
\end{algorithm}

\begin{fact}
The strategy of restarts in Algorithm \ref{a:RSV-1} preserves the
summability of the overall series of step-sizes. This is so because
even if during the iterative process, the largest step-sizes allowed
in each of the sets \textup{$W_{r}$} were taken, then the infinite
sequence of all step-sizes
\begin{equation}
\left(\left(\alpha^{r+\ell}\right)_{\ell=0}^{W_{r}-1}\right)_{r=0}^{\infty},\label{e:bound_seq_restarts}
\end{equation}
forms a bounded series. We have, since $\alpha\in{]0,1[}$, 
\begin{equation}
\sum_{r=0}^{\infty}\sum_{\ell=0}^{W_{r}-1}\alpha^{r+\ell}\leq\sum_{r=0}^{\infty}\alpha^{r}\sum_{\ell=0}^{\infty}\alpha^{\ell}=\frac{1}{(1-\alpha)^{2}}.\label{e:restbound}
\end{equation}
Hence, since only step-sizes leading to expected superior values of
the target function are allowed by~\eqref{e:supphi} (line 8 of Algorithm\ \ref{a:RSV-1}),
each of the step-sizes taken will be smaller than the corresponding
one in the sequence~\eqref{e:bound_seq_restarts}, so its sum will
always be smaller than $1/(1-\alpha)^{2}$ and will, thus, define
bounded perturbations.
\end{fact}

\smallskip{}

For some applications, the SM with restarts is very useful, notably
outperforming the current SM without restarts (see Section \ref{exp:3}).

\section{Second development: A superiorized algorithm for subvectors in the
split minimization problem\label{Sect:3}}

We develop here a superiorized algorithm for tackling the data of
the SMP in Problem~\ref{prob:smp-1} when $C:=\bigcap_{s=1}^{p}C_{s}$
and $Q:=\bigcap_{t=1}^{q}Q_{t},$ where $p$ and $q$ are two integers
and $\{C_{s}\}_{s=1}^{p}$ and $\{Q_{t}\}_{t=1}^{q}$ are two families
of closed and convex sets with nonempty intersections in $\mathbb{R}^{n}$
and $\mathbb{R}^{m}$, respectively. To ease the discussion we will
refer  here to $\mathbb{R}^{n}$ and $\mathbb{R}^{m},$ as the ``$x$-space''
and the ``$y$-space'', respectively.

\subsection{The SMP with subvectors}

In some situations of practical interest, the minimization problem
in (\ref{e:svip2}) should be independently applied to subvectors
of the $y$-space. We discuss an instance in the field of radiation
therapy treatment planning where this is significant in Section \ref{Sect:4}
below.

For simplicity and without loss of generality, we assume that the
subvectors are in consecutive order. The $m\times n$ real matrix
$A$ is divided into $B$ blocks and is represented by
\begin{equation}
A:=\left(\begin{array}{c}
A_{1}\\
A_{2}\\
\vdots\\
A_{B}
\end{array}\right)\label{eq:blockA}
\end{equation}
 where, for each $b\in B,$ $A_{b}\in\mathbb{R}^{m_{b}\times n}$
are blocks of rows of the matrix $A$, with $\sum_{b=1}^{B}m_{b}=m$.
Thus, any vector $y:=Ax$ is of the form
\begin{equation}
y=\left(\begin{array}{c}
y^{1}\\
y^{2}\\
\vdots\\
y^{B}
\end{array}\right)=\left(\begin{array}{c}
A_{1}\\
A_{2}\\
\vdots\\
A_{B}
\end{array}\right)x
\end{equation}
where $y^{b}\in\mathbb{R}^{m_{b}}$ are subvectors of $y\in\mathbb{R}^{m}.$
\begin{problem}
\label{prob:subvectors}\textbf{The}\textbf{\emph{ }}\textbf{SMP with
subvectors}. Given two families of closed and convex sets $\{C_{s}\}_{s=1}^{p}\subseteq\mathbb{R}^{n}$
and $\{Q_{t}\}_{t=1}^{q}\subseteq\mathbb{R}^{m}$ such that $C:=\bigcap_{s=1}^{p}C_{s}\neq\emptyset$
and $Q:=\bigcap_{t=1}^{q}Q_{t}\neq\emptyset$, an $m\times n$ real
matrix $A$ in the form (\ref{eq:blockA}) for a given integer $B$,
a convex function $f:\mathbb{R}^{n}\to\mathbb{R}$ and convex functions
and $\phi_{b}:\mathbb{R}^{m_{b}}\to\mathbb{R}$, for $b=1,2,\ldots,B$,
find
\begin{equation}
\begin{aligned}x\in C\;\text{such that }x^{*}\in\argmin\{f(x)\mid x\in C\},\;\text{ and such that}\;\end{aligned}
\label{e:svipsub3}
\end{equation}
\begin{equation}
\begin{aligned}y^{*}:=Ax^{*}\in Q\;\text{and}\;y^{b*}\in\argmin\{\phi_{b}(y^{b})\mid y\in Q\}\;\text{ for all }\;b\in{\{1,2,\ldots,B\}}.\end{aligned}
\label{e:svipsub4-1}
\end{equation}
\end{problem}

Our algorithm, presented below, can also handle subvectors in the
$x$-space, but for simplicity we restrict ourselves here to subvectors
in the $y$-space.

\subsection{Reformulation in the product space}

To work out a superiorization method for the\textbf{ }data of the
SMP with subvectors in Problem \ref{prob:subvectors} we look at a
multiple sets split feasibility problem (MSSFP), see, e.g., \cite{masad-2007},
as follows.
\begin{problem}
\label{prob:mssfp}\textbf{The multiple sets split feasibility problem
(MSSFP)}. Given $C:=\bigcap_{s=1}^{p}C_{s}$and $Q:=\bigcap_{t=1}^{q}Q_{t},$
where $p$ and $q$ are two integers, and $\{C_{s}\}_{s=1}^{p}$ and
$\{Q_{t}\}_{t=1}^{q}$ are two families of closed and convex sets
with nonempty intersections each in $\mathbb{R}^{n}$ and $\mathbb{R}^{m}$,
respectively, and an $m\times n$ real matrix $A$, find 
\begin{equation}
x^{*}\in C=\bigcap_{s=1}^{p}C_{s}\;\text{such that\;}y^{*}:=Ax^{*}\in Q=\bigcap_{t=1}^{q}Q_{t}.
\end{equation}
\end{problem}

This is a generalization of the \textit{split feasibility problem}
(SFP) that occurs when $p=q=1$ in the MSSFP. The SFP, which plays
an important role in signal processing, in medical image reconstruction
and in many other applications, was introduced by Censor and Elfving
\cite{CE-SFP-1994} in order to model certain inverse problems. Since
then, many iterative algorithms for solving the SFP have been proposed
and analyzed. See, for instance, the references given in \cite{Reich-split-2021}
or consult the section ``A brief review of 'split problems' formulations
and solution methods'' in \cite{brooke-2020}.

Our proposed algorithm deals with an equivalent reformulation of Problem~\ref{prob:mssfp}
in the product space $\mathbb{R}^{n}$$\times$$\mathbb{R}^{m}$.
Adopting the notation that quantities in the product space are denoted
by boldface symbols, we define the sets 
\begin{equation}
\mathbf{C}:=\left(\bigcap_{s=1}^{p}C_{s}\right)\times\left(\bigcap_{t=1}^{q}Q_{t}\right)\quad\text{and}\quad\mathbf{V}:=\lbrace\mathbf{z}=(x,y)\in\mathbb{R}^{n}\times\mathbb{R}^{m}\mid Ax=y\rbrace.
\end{equation}
Note that the projection onto $\mathbf{V}$ is given by the expression $P_{\mathbf{{V}}}= I - Z^{T}(ZZ^{T})^{-1}Z$,
with  $Z:=[A,-I]$ and where  $I$ stands for  the  identity matrix of appropriate order.
Then, Problem \ref{prob:mssfp} is equivalent to the problem:
\begin{equation}
\text{find a point}\;\mathbf{z}^{*}\in\mathbf{C}\cap\mathbf{V}.\label{e:fp_ps}
\end{equation}
Without loss of generality, we assume that $p=q,$ since otherwise
the whole space (or one particular set) could be added repeatedly
as a constraint until both indices are equal. Since the projection
of a Cartesian product is the Cartesian product of the projections~\cite[Lemma~1.1]{pierra},
the following implementation of the method of alternating projections
can be employed to solve~\eqref{e:fp_ps}. We consider this as our
basic algorithm for the superiorization method for subvectors.

\smallskip{}

\begin{algorithm}
\textbf{Initialization:}~Choose an arbitrary initialization point
$x^{0}\in\mathbb{R}^{n}$. Set $y^{0}=Ax^{0}$\; \textbf{Iterative
Step:}~Given the current iterate $\mathbf{z}^{k}=(x^{k},y^{k})$,
calculate the next iterate $\mathbf{z}^{k+1}$ by 
\begin{equation}
\mathbf{z}^{k+1}=P_{\mathbf{V}}\left(\left(P_{C_{p}}\times P_{Q_{p}}\right)\cdots\left(P_{C_{1}}\times P_{Q_{1}}\right)(\mathbf{z}^{k})\right).\label{e:basic_svip}
\end{equation}
\caption{Basic Algorithm for~the MSSFP.}
\label{a:BA_svip}
\end{algorithm}

\smallskip{}

In order to construct a superiorized version of Algorithm~\ref{a:BA_svip}
that can cope with the\textbf{ }data of the SMP with subvectors in
Problem \ref{prob:subvectors}, we need to establish at each iteration
some appropriate perturbations that will steer the algorithm to a
superiorized solution. For this, we note that the vector $y^{k}$
inside $\mathbf{z}^{k}$ in~\eqref{e:basic_svip} is expressed as
\begin{equation}
y^{k}=\left(\begin{array}{c}
y^{1,k}\\
y^{2,k}\\
\vdots\\
y^{B,k}
\end{array}\right)=\left(\begin{array}{c}
A_{1}x^{k}\\
A_{2}x^{k}\\
\vdots\\
A_{B}x^{k}
\end{array}\right).
\end{equation}
Thus, we declare our perturbation vector to be 
\begin{equation}
\left(\begin{array}{c}
x^{k}\\
y^{k}
\end{array}\right)+\eta_{k}\left(\begin{array}{c}
u^{k}\\
v^{k}
\end{array}\right),\quad\text{for all}\;k=0,1,2,\ldots,\label{e:svip_pert}
\end{equation}
where $\{\eta_{k}\}_{k=0}^{\infty}$ is a nonnegative summable sequence,
$u^{k}$ is a nonascending vector for $f$ at $x^{k}$ and
\begin{equation}
v^{k}=\left(\begin{array}{c}
v^{1,k}\\
v^{2,k}\\
\vdots\\
v^{B,k}
\end{array}\right),
\end{equation}
with each $v^{b,k}$ being a nonascending vector for $\phi_{b}$ at
the point $A_{b}x^{k}$, for all $b\in{\{1,2,\ldots,B\}}$. The complete
pseudo-code of the superiorized version of the basic Algorithm~\ref{a:BA_svip}
with perturbations of the form given by~\eqref{e:svip_pert} is shown
in Algorithm~\ref{a:SV_SPM}.

\begin{algorithm}
\textbf{Initialization:}~Choose $x^{0}\in\mathbb{R}^{n}$, a summable
nonnegative sequence $\{\eta_{\ell}\}_{\ell=0}^{\infty}$ and a positive
integer $N$\; \Set{$y^{b,0}=A_{b}x^{0}$ for $b\in{\{1,2\text{,}\ldots,B\}}$};
\Set{$k=0$ and $\ell=-1$}\; \Repeat{ \Set{$x^{k,0}=x^{k}$}\;
\Set{$y^{b,k,0}=y^{b,k}$ for $b\in{\{1,2\text{,}\ldots,B\}}$}\;
\For{$j=0$ \KwTo $N-1$}{ \Set{$v_{x}^{k,j}$ to be a nonascending
vector for $f$ at $x^{k,j}$}\; \Set{$v_{y}^{b,k,j}$ to be a
nonascending vector for $\phi_{b}$ at the point $y^{b,k,j}$ for
all $b$ }\; \Set{$\ell=\ell+1$}\; \While{$f(x^{k,j}+\eta_{\ell}v_{x}^{k,j})>f(x^{k})$
\textbf{or there exists $b$ with} $\phi_{b}(y^{b,k,j}+\eta_{\ell}v_{y}^{b,k,j})>\phi_{b}(y^{b,k})$}{
\Set{$\ell=\ell+1$}\; } \Set{$x^{k,j+1}=x^{k,j}+\eta_{\ell}v_{x}^{k,j}$}\;
\Set{$y^{b,k,j+1}=y^{b,k,j}+\eta_{\ell}v_{y}^{b,k,j}$}\; } \Set{$(x^{k+1},y^{1,k+1},\ldots,y^{B,k+1})=P_{\mathbf{V}}\left(\left(P_{C_{p}}\times P_{Q_{p}}\right)\cdots\left(P_{C_{1}}\times P_{Q_{1}}\right)\left(x^{k,N},y^{1,k,N},\ldots,y^{B,k,N}\right)\right)$
}\; \Set{$k=k+1$}\; } \caption{Superiorized Algorithm for the data of the SMP with subvectors}
\label{a:SV_SPM}
\end{algorithm}

 Since the method of alternating projections is bounded
perturbation resilient~\cite{sm2007}, Algorithm \ref{a:SV_SPM}
will converge to a solution of the feasibility problem~\eqref{e:fp_ps}.
Moreover, by the nature of the SM, the algorithm is expected to converge
to a point $\mathbf{z}^{*}=(x^{*},y^{*})$ which will be superior
with respect to $f$ for the component $x$ in the $x$-space, and
with respect to $\phi_{b}$ for the $b$-th subvector in the $y$-space,
for $b=1,2,\ldots,B$.

\section{Numerical experiments\label{Sect:4}}

Our aim in this paper is not to compare the superiorization
method with constrained optimization methods, moreover, the SM is
not a method intended to solve exact constrained optimization problems.
Such comparisons were done elsewhere, see, e.g., \cite{WHare} or
\cite{supvssubg}. Our goal is to show how the SM can be improved
and this is achieved by comparing the SM with and without restarts
and with and without perturbations. We present
our results of numerical experiments performed on three different
problems. The first two problems are simple illustrative examples.
The computational performance of the SM with restart algorithms, proposed
here, can be substantiated with exhaustive testing of the possible
specific variants permitted by the general framework and their various
user-chosen parameters. This should be done on larger problems, preferably
within the context of a significant real-world application. Therefore,
our third problem addresses an actual situation arising in the real-world
application of intensity-modulated radiation therapy (IMRT) treatment
planning.

The purpose of the first example is to illustrate the potential benefits
of superiorization with restarts for finding a point with reduced
norm in the intersection of two convex sets. We first consider the
case of two balls and then explore the case of two half-spaces presented
in Remark~\ref{rem:failure}, in which superiorization did not achieve
its purpose for a particular setting.

In the second example we illustrate the behavior of Algorithm \ref{a:SV_SPM}
in a simple setting with $C,Q\subset\mathbb{R}^{2}$, where each of
the sets is an intersections of three half-spaces.

Finally, the last experiment demonstrates the benefits of superiorization
with restarts in a difficult realistic setting in IMRT, where a large-scale
multiobjective optimization problem arises. All tests were run on
a desktop of Intel Core i7-4770 CPU 3.40GHz with 32GB RAM, under Windows
10 (64-bit).

\subsection{The benefits of superiorization with restarts}

\label{exp:1}

Consider the problem of finding the minimum norm point in the intersection
of two balls $A$ and $B$ in the Euclidean 2-dimensional space, so
$\phi:\mathbb{R}^{2}\to\mathbb{R}$ is given by $\phi(x):=\frac{1}{2}\|x\|^{2}$
for $x\in\mathbb{R}^{2}$. The underlying feasibility problem can
be solved by the method of alternating projections, which we chose
as the basic algorithm. Hence, the feasibility-seeking algorithmic
operator used in our computations is 
\begin{equation}
\mathbf{T}=P_{B}P_{A},
\end{equation}
where $P_{A}$ and $P_{B}$ denote the projection operators onto the
balls $A$ and $B,$ respectively. We tested the method of alternating
projections, its superiorized version (with two different kernels
$\alpha=0.6$ and $\alpha=0.999$) and its superiorized version with
restarts (with $\alpha=0.6$ and $W_{r}=50$, for all $r\geq0$).
We set $N=1$ in all the superiorized algorithms. The nonascending
directions were taken as $v^{k}:=-\frac{\nabla\phi(y^{k})}{\|\nabla\phi(y^{k})\|}=-\frac{y^{k}}{\|y^{k}\|}$.\smallskip{}

\begin{figure}[H]
\includegraphics[width=1\textwidth]{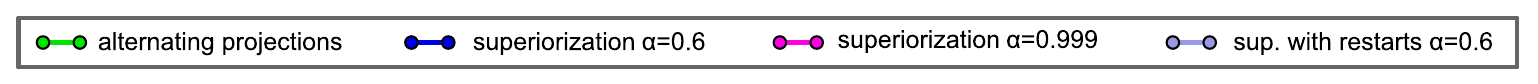}

\includegraphics[width=0.49\textwidth]{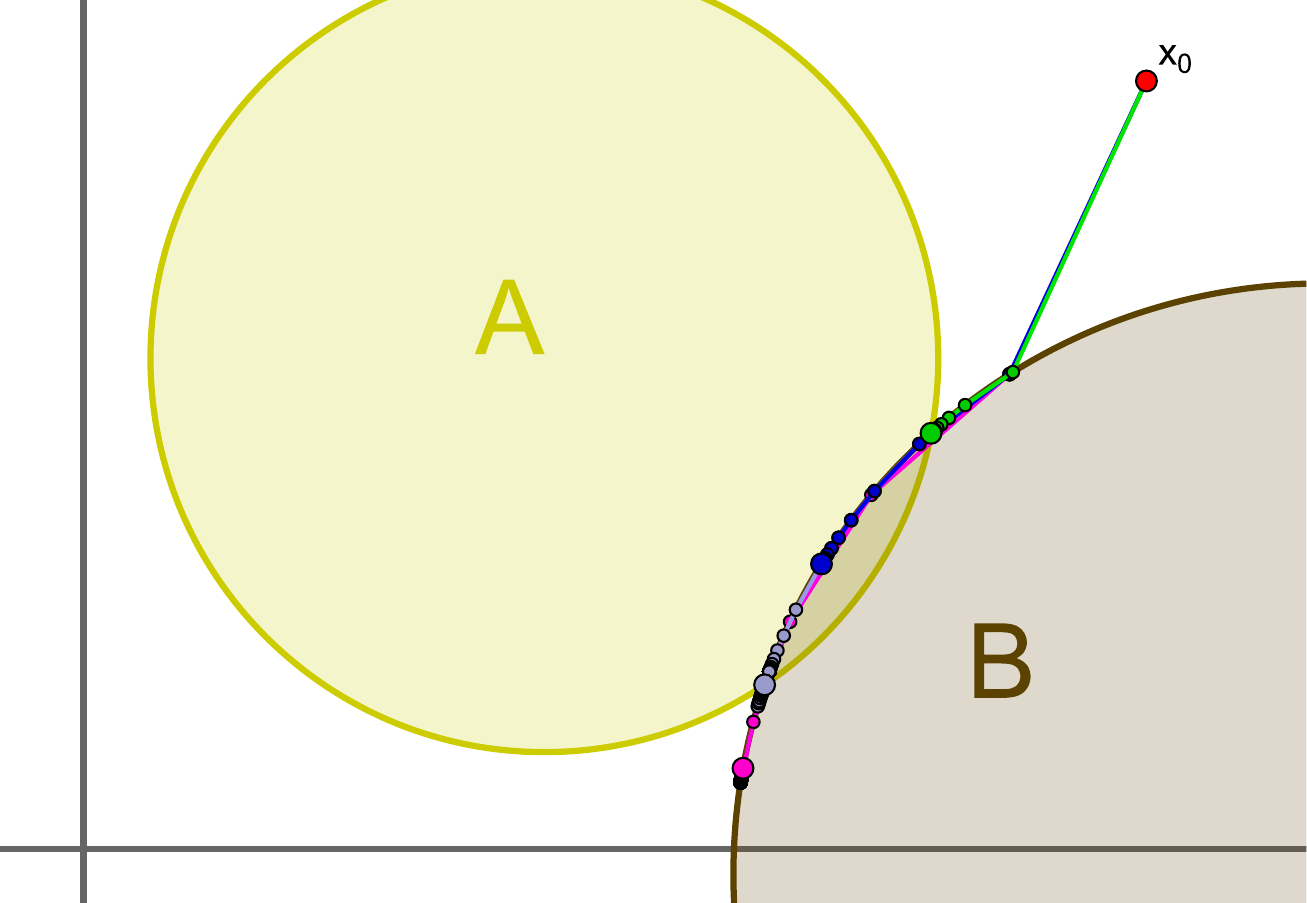}
\includegraphics[width=0.49\textwidth]{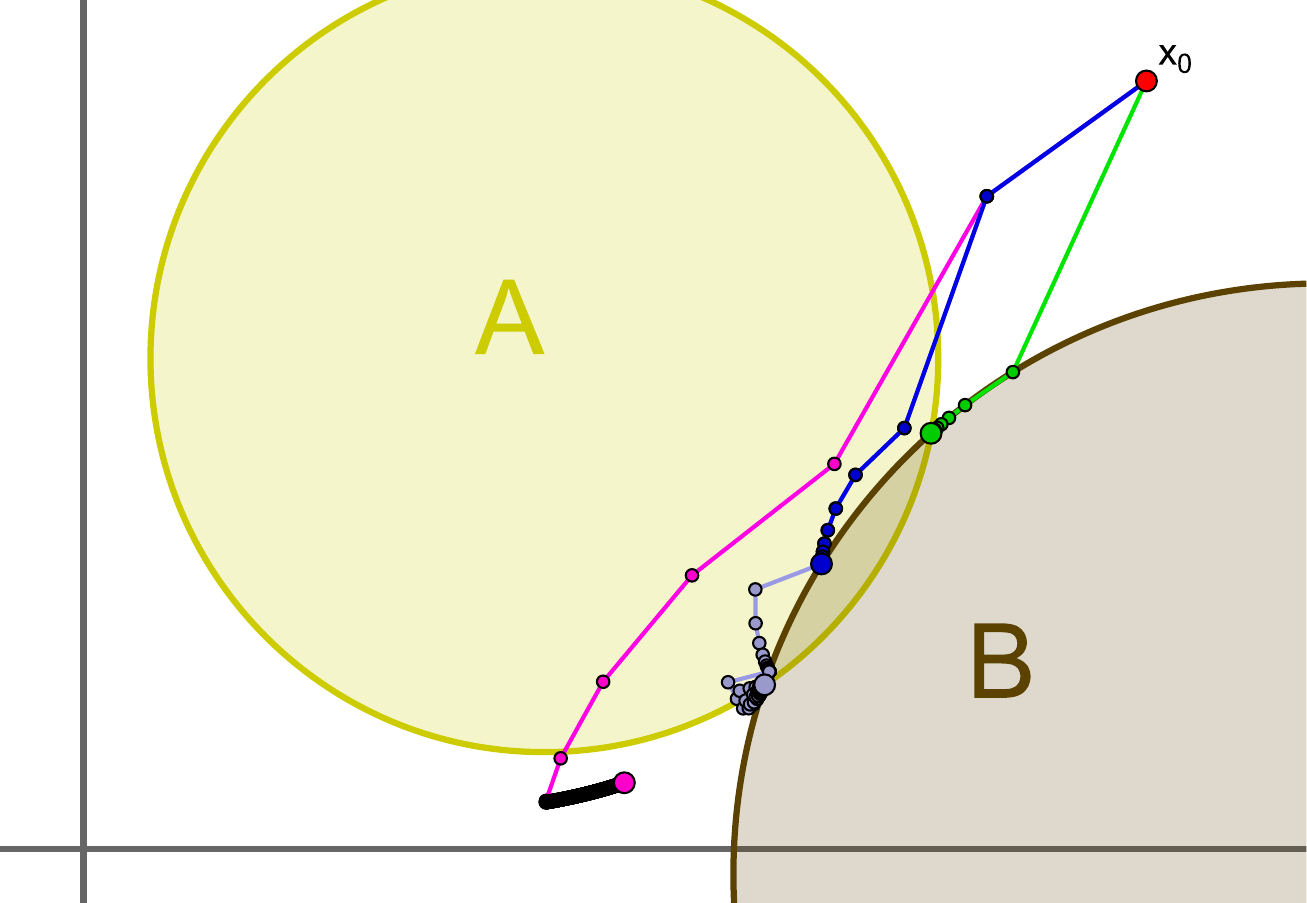}
\caption{Behavior of the different algorithms considered applied to the data
of the problem of finding the minimum norm point in the intersection
of the balls $A$ and $B$. The figures show the first 500 points
in the sequence of iterates $\{y^{k}\}_{k=0}^{500}$ (left) and the
sequence of perturbed iterations before applying the algorithmic operator
$\{y^{k,N}\}_{k=0}^{500}$ (right) of each algorithm.}
\label{fig:1}
\end{figure}

The behavior of these algorithms is shown in Figure~\ref{fig:1}.
On the left we represent 500 iterations generated by each algorithm.
On the right, we plot the sequence of perturbations obtained before
applying the algorithmic operator, that is, we draw the points $\left\{ y^{k,N}\right\} _{k=0}^{500}$.
This sequence coincides with the sequence of iterates $\left\{ y^{k}\right\} _{k=0}^{500}$
in the case when no perturbations are performed at all and only the
basic algorithm works, while it coincides with the sequence $\left\{ y^{k}+\eta_{k}v^{k}\right\} _{k=0}^{500}$
for the superiorized algorithms.

As expected, the method of alternating projections converges to a
point in the intersection which is not desirable according to the
task of reducing the target function value (the squared norm). Superiorization
with kernel $\alpha=0.6$ reaches a better point than the output of
the basic algorithm, but is yet far from the solution to the problem.
This might well be due to the step-sizes not being big enough for
the perturbations to steer the algorithm to a proper function reduction.

Taking $\alpha=0.999$ in the standard superiorized version of the
algorithm results in a very slow convergence of the algorithm, as
can be observed on the right figure in Figure~\ref{fig:1}. These
deficiencies are resolved by considering superiorization with restarts,
which achieves fast convergence to a solution with reduced norm in
the intersection. 

The example presented in Section~\ref{rem:failure} is artificial
in the sense that the vectors defining the half-spaces are orthogonal
and the starting point was chosen in a particular region of the plane
which was less favorable to the superiorized algorithm. Also, the
value of the kernel was chosen to be small ($\alpha=0.5$), to aggravate
the vanishing effect of the perturbations.

To investigate what happens with random data, we ran an experiment
generating $1\,000\,000$ pairs of half-spaces of the form $A:=\left\{ x\in\mathbb{R}^{2}\mid\left\langle c_{A},x\right\rangle \leq b_{A}\right\} $
and $B:=\left\{ x\in\mathbb{R}^{2}\mid\left\langle c_{B},x\right\rangle \leq b_{B}\right\} $
where the vectors $c_{A},c_{B}\in\mathbb{R}^{2}$ were randomly chosen
and then normalized, and $b_{A},b_{B}\in(-1,0)$ (to ensure that $(0,0)^{T}\not\in A\cap B$).
For each pair of half-spaces, we generated a random starting point
$y^{0}\in[-1,1]^{2}$ such that $y^{0}\not\in A\cap B$.

Then, we ran from $y^{0}$ the basic algorithm (feasibility-seeking
alternating projections), its superiorized version with kernel $\alpha\in\left\{ 0.5,0.6,0.7,0.8,0.9\right\} $
and its superiorized version with restarts with $W_{r}=20$ and $N=1$.
The results are summarized in Table~\ref{tab:halfspaces}. In this
table AP stands for ``alternating projections'', Sup stands for
``the superiorized version'', and Sup. Res. stands for ``the superiorized
version with restarts''. We count that one method is better than
the other when the norm of its solution is smaller than the norm of
the second method's output minus $10^{-3}$. With kernel $\alpha=0.5$,
the superiorized algorithm failed to obtain a solution with lower
norm than the basic algorithm in $12\,931$ of the $1\,000\,000$ instances;
with kernel $\alpha=0.9$, this number was reduced to 122. The superiorized
algorithm with restarts with kernel $\alpha=0.7$ only failed to get
a better solution than the basic algorithm in 21 instances. Remarkably,
when the kernels $\alpha\in\left\{ 0.8,0.9\right\} $ were used, superiorization
with restarts always reached the same or a better solution than both
the basic and the superiorized algorithms without restarts.

\begin{table}[H]
\begin{centering}
\begin{tabular}{|c|c|c|c|r@{\extracolsep{0pt}.}l|c|r@{\extracolsep{0pt}.}l|}
\cline{2-9} \cline{3-9} \cline{4-9} \cline{5-9} \cline{7-9} \cline{8-9} 
\multicolumn{1}{c|}{} & \multicolumn{2}{c|}{AP vs Sup.} & \multicolumn{3}{c|}{AP vs Sup. Res.} & \multicolumn{3}{c|}{Sup. vs Sup. Res.}\tabularnewline
\hline 
Kernel & \multirow{1}{*}{AP} & \multicolumn{1}{c|}{Sup.} & AP & Sup& Res. & Sup. & Sup& Res.\tabularnewline
\hline 
$\alpha=0.5$ & 1.29\% & 56.17\% & 0.08\% & 57&2\% & 0.01\% & 16&9\%\tabularnewline
\hline 
$\alpha=0.6$ & 0.73\% & 56.63\% & 0.02\% & 57&26\% & 0.001\% & 10&78\%\tabularnewline
\hline 
$\alpha=0.7$ & 0.32\% & 56.96\% & 0.002\% & 57&28\% & 0\% & 6&1\%\tabularnewline
\hline 
$\alpha=0.8$ & 0.10\% & 57.17\% & 0\% & 57&29\% & 0\% & 2&86\%\tabularnewline
\hline 
$\alpha=0.9$ & 0.01\% & 57.27\% & 0\% & 57&3\% & 0\% & 0&68\%\tabularnewline
\hline 
\end{tabular}
\par\end{centering}
\caption{For each pair-wise comparative of methods and kernel choice, the numbers
inside the table are the percentage of the $1\,000\,000$ runs in which
each method obtains a solution with lower norm than the other one
with which it is compared.\label{tab:halfspaces}}
\end{table}

\subsection{Behavior of the superiorized algorithm for the data of the SMP}

In this section we present another illustrative example of the performance
of Algorithm~\ref{a:SV_SPM}. To be able to display the iterates,
we let both the $x$-space and the $y$-space in Problem~\eqref{e:svipsub3}-\eqref{e:svipsub4-1}
be the Euclidean 2-dimensional spaces. We take $C$ as the intersection
of the three half-spaces given by $C_{1}:=\left\{ (x_{1},x_{2})\in\mathbb{R}^{2}\mid x_{1}+x_{2}\leq10\right\} $,
$C_{2}:=\left\{ (x_{1},x_{2})\in\mathbb{R}^{2}\mid-13x_{1}+3x_{2}\leq-26\right\} $
and $C_{3}:=\left\{ (x_{1},x_{2})\in\mathbb{R}^{2}\mid x_{2}\geq1\right\} $.

We let $A$ be the rotation matrix by an angle of $\pi/2$ and $Q$
be the image of $C$ under $A$ (i.e., the intersection of the half-spaces
$Q_{1}:=\left\{ (y_{1},y_{2})\in\mathbb{R}^{2}\mid y_{1}-y_{2}\leq-10\right\} $,
$Q_{2}:=\left\{ (y_{1},y_{2})\in\mathbb{R}^{2}\mid-3y_{1}-13y_{2}\leq-26\right\} $
and $Q_{3}:=\left\{ (y_{1},y_{2})\in\mathbb{R}^{2}\mid y_{1}\leq-1\right\} $).
The function to be reduced in the $x$-space is the value of the second
component $f(x_{1},x_{2}):=x_{2}$, whereas in the $y$-space, we
aim to find a point with increased first and second components (that
is, $B:=2$, $\phi_{1}(y_{1}):=-y_{1}$ and $\phi_{2}(y_{2}):=-y_{2}$).

In other words, we want to tackle with the SM the data of the split
minimization problem given by 
\begin{equation}
\begin{aligned}\text{find }x^{*}=\left(\begin{array}{c}
x_{1}^{*}\\
x_{2}^{*}
\end{array}\right)\in C\text{ such that }x^{*}\in\argmin{\left\{ x_{2}\;\Bigl|\;\left(\begin{array}{c}
x_{1}\\
x_{2}
\end{array}\right)\in C\right\} }\end{aligned}
\label{eq:exp2-1}
\end{equation}

\begin{equation}
\begin{aligned}\textup{and such that }\left(\begin{array}{c}
-x_{2}^{*}\\
x_{1}^{*}
\end{array}\right)\in Q,\quad-x_{2}^{*}\geq y_{1}\text{ and }x_{1}^{*}\geq y_{2}\text{ for all }\left(\begin{array}{c}
y_{1}\\
y_{2}
\end{array}\right)\in Q,\end{aligned}
\label{eq:exp2-2}
\end{equation}
with $C:=\bigcap_{i=1}^{3}C_{i}$ and $Q:=\bigcap_{j=1}^{3}Q_{j.}$

\begin{figure}[ht!]
\centering \includegraphics[width=0.45\textwidth]{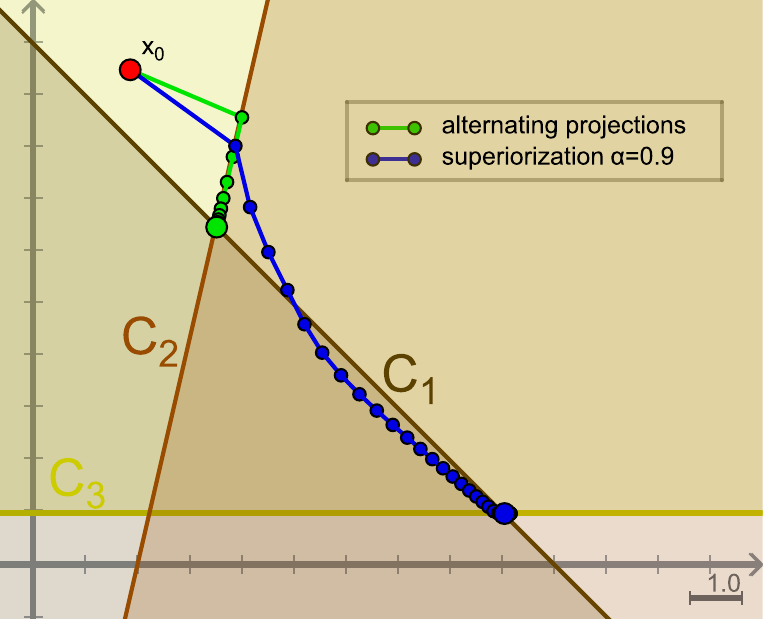}
\hspace{0.01\textwidth} \includegraphics[width=0.45\textwidth]{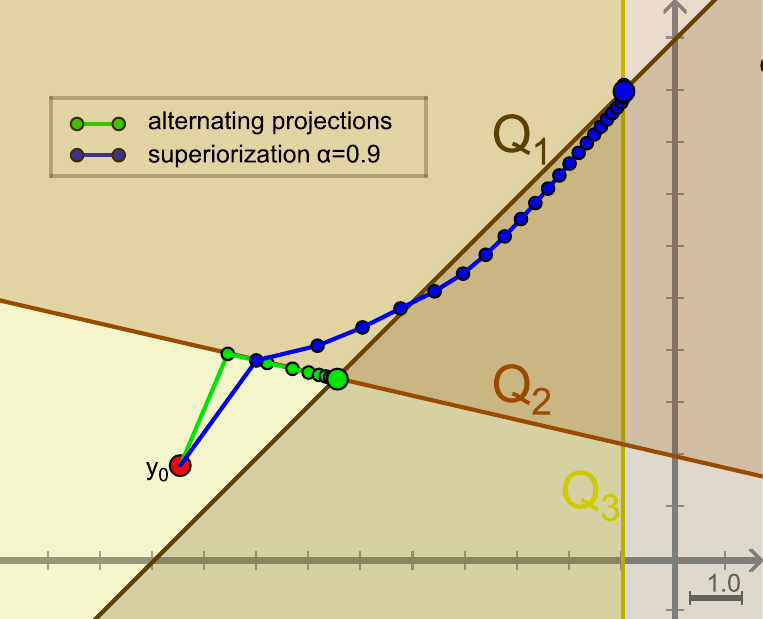}
\caption{Performance of the method of alternating projections and its superiorized
version applied in the setting of problem \eqref{eq:exp2-1}-\eqref{eq:exp2-2}.
The left image displays the iterates in the $x$-space, while the
right one shows the iterates in the $y$-space.}
\label{fig:2}
\end{figure}

By looking at Figure~\ref{fig:2}, one easily identifies that the
point $(9,1)^{T}$, obtained as the intersection of the lines $x_{1}+x_{2}=10$
and $x_{2}=1$, is the unique solution to the SMP with the above described
data.

Again, for the SM we choose the method of alternating projections
as the basic algorithm. Consequently, the algorithmic operator that
we use is 
\begin{equation}
\mathbf{T}:=P_{\textbf{V}}\circ\left(P_{C_{3}}\times P_{Q_{3}}\right)\circ\left(P_{C_{2}}\times P_{Q_{2}}\right)\circ\left(P_{C_{1}}\times P_{Q_{1}}\right),
\end{equation}
where we recall that $\mathbf{V}:=\left\{ (x,y)\in\mathbb{R}^{2}\times\mathbb{R}^{2}\mid Ax=y\right\} $.
In our experiment, we performed 50 iterations of both the basic algorithm
and its superiorized version, taking $\alpha=0.9$ as the kernel for
generating the step-sizes of the perturbations and $N=1$. The nonascending
vectors were taken as $v_{x}^{k}:=(0,1)^{T}$ for the perturbations
in the $x$-space, and $v_{y}^{1,k}:=1$ and $v_{y}^{2,k}:=1$ in
the $y$-space.

Figure~\ref{fig:2} shows that, while the method of alternating projection
converges to the closest point to the starting point in the intersection
in each of the spaces, the superiorized algorithm converges to the
solution of the SMP.

\subsection{A demonstrative example in IMRT}

\label{exp:3}

In this section, we test our SM with restarts algorithm in a sophisticated
multiobjective setting motivated from a split minimization problem
in the field of intensity-modulated radiation therapy (IMRT) treatment
planning. IMRT is a radiation therapy that manipulates particle beams
(protons or photons or others) of varying directions and intensities
that are directed toward a human patient to achieve a goal of eradicating
tumorous tissues, henceforth called ``tumor structures'', while
keeping healthy tissues, called ``organs-at-risk'' below certain
thresholds of absorbed dose of radiation. The beams are projected
onto the \emph{region of interest} from different angles.

Many review papers in this field are available, see, e.g., \cite{cho-marks-2000,unkelbach-2012,imrt-review-2018,nonconv-2021,TRPMS-2022}
and references therein.

\subsubsection{The fully-discretized model of the inverse problem of IMRT}

In the fully-discretized model of the inverse problem of IMRT each
external radiation beam is discretized into a finite number of ``beamlets''
(also called ``pencil-beams'' or ``rays'') along which the particles
(i.e., their energies) are transmitted. Let all beamlets from all
directions be indexed by $i=1,2,\ldots,n,$ and denote the ``intensity''
irradiated along the $i$-th beamlet by $x_{i}$. The vector $x=(x_{i})_{i=1}^{n}\in\mathbb{R}^{n}$
is called the ``intensities vector''.

The 2-dimensional (2D) cross-section\footnote{Everything presented here can easily be extended to 3D wherein the
pixels are replaced by voxels. The choice of the 2D case just makes
the presentation simpler.} of the irradiated body is discretized. Assume that the cross-section
is covered by a square that is discretized into a finite number of
square pixels. This creates an $M\times M$ array of pixels. Let all
pixels be indexed by $j=1,2,\ldots,m,$ with $m=M^{2}$ and let $y_{j}$
denote the ``dose'' of radiation absorbed in the $j$-th pixel.
The vector $y=(y_{j})_{j=1}^{m}\in\mathbb{R}^{m}$ is called the ``dose
vector''.

The ``intensities space'' $\mathbb{R}^{n}$ and the ``dose space''
$\mathbb{R}^{m},$ defined above are the ``$x$-space'' and the
``$y$-space'', respectively, mentioned at the beginning of Section~\ref{Sect:3}.
The physics of the model assumes that there exists an $m\times n$
real matrix $A=(a_{ij})_{i=1,j=1}^{n,m}$(sometimes called ``the
dose matrix'') through which the intensities of the beamlets and
the absorbed doses in pixels are related via the equation
\begin{equation}
Ax=y.
\end{equation}

Each element $a_{ij}$ in $A$ is the dose absorbed in pixel $j$
due to a unit of intensity along the $i$-th beamlet. This means that
\begin{equation}
\sum_{i=1}^{n}a_{ij}x_{i}
\end{equation}
is the total dose absorbed in pixel $j$ due to an intensity vector
$x$. With these notions in mind we consider the following feasibility-seeking
problem of the fully-discretized inverse problem of IMRT.
\begin{problem}
\label{prob:IMRT0}\textbf{The feasibility-seeking problem of the
fully-discretized inverse problem of IMRT}. Let $\mathbb{R}^{n}$
and $\mathbb{R}^{m}$ be the ``intensities space'' and the ``dose
space'' (henceforth called the ``$x$-space'' and the ``$y$-space'')
respectively. Let $A=(a_{ij})_{i=1,j=1}^{n,m}$ be the dose matrix
mapping the $x$-space onto the $y$-space. For $\ell=1,2,\ldots,L$,
denote by $T_{\ell}\subseteq\{1,2,\ldots,m\}$ the set of pixels corresponding
to the $\ell$-th tumor structure in the region of interest. For $r=1,2,\ldots,R$,
denote by $S_{r}\subseteq\{1,2,\ldots,m\}$ the set of pixels corresponding
to the $r$-th organ-at-risk. Set $\underline{e}$ and $\overline{e}$
the lower and upper bounds for the available beamlets intensities.
Let $\underline{d}_{r}$ and $\overline{d}_{r}$, and $\underline{c}_{\ell}$
and $\overline{c}_{\ell}$ be the lower and upper bounds for
the dose deposited in each pixel of the $r$-th organ at risk and
of the $\ell$-th tumor,
respectively.

The task is to find an intensities vector $x$ such that 
\begin{equation}
\begin{aligned}\underline{c}_{\ell}\leq\sum_{i=1}^{n}a_{ij}x_{i}\leq\overline{c}_{\ell}, & \quad\text{ for all }j\in T_{\ell},\quad\ell\in{\{1,2,\ldots,L\}},\\
\underline{d}_{r}\leq\sum_{i=1}^{n}a_{ij}x_{i}\leq\overline{d}_{r}, & \quad\text{ for all }j\in S_{r},\quad r\in{\{1,2,\ldots,R\}},\\
\underline{e}\leq x_{i}\leq\overline{e}, & \quad\text{ for all }i\in{\{1,2,\ldots,n\}}.
\end{aligned}
\end{equation}
\end{problem}

Problem~\ref{prob:IMRT0} is a linear feasibility problem. Usually,
$\underline{e}=\underline{d}_{r}=0$ and the significant bounds are
$\underline{c}_{\ell}$ for tumor structures and $\overline{d}_{r}$
for organs-at-risk. A pair $(x^{*},y^{*})$ such that $x^{*}$ is
a solution of Problem \ref{prob:IMRT0} and $y^{*}=Ax^{*}$ will be
henceforth called ``a treatment plan'' for the IMRT inverse planning
problem.

\subsubsection{The quest for smoothness and uniformness}

In the IMRT inverse planning problem there is an advantage to generating
treatment plans with intensity vectors $x^{*}$ whose subvectors,
related to parallel beamlets from the same beam, will be as ``smooth''
as possible and with dose vectors $y^{*}$ whose subvectors, related
to specific organs (a.k.a. ``structures''), that will be as ``uniform''
as possible.

For the intensity vectors $x^{*},$ ``smoothness'' of subvectors,
related to parallel beamlets from the same beam, means that the real
numbers that are the individual intensities $x_{i}^{*}$, in each
subvector separately, would be as close to each other as possible,
subject to the constraints of Problem \ref{prob:IMRT0}. Such smoothness
will allow for less extreme movements of the ``multileaf collimator''
(MLC)\footnote{A multileaf collimator is a beam-limiting device that is made of individual
``leaves'' of a high atomic numbered material, usually tungsten,
that can move independently in and out of the path of a radiotherapy
beam in order to shape (i.e., modulate) it and vary its intensity.
See, e.g., \cite{MLC}.} that modulates the parallel beamlets from the same beam.

For the dose vectors $y^{*},$ ``uniformness'' of subvectors, related
to specific organs means that the real numbers that are the individual
doses $y_{j}^{*}$, in each pixel of the subvector would be as close
to each other as possible, subject to the constraints of Problem \ref{prob:IMRT0}.
Such uniformness will guarantee uniformness of the dose deposited
within each organ separately and help to avoid the presence of hot-
and cold-spots in the dose distribution in each organ. See, e.g.,
\cite{cold-hot}.

Each of these aims can be achieved by attempting to minimize or just
reduce the \emph{total variation }(TV)\emph{-norm} of the associated
subvectors, see, e.g., \cite{TV-2010} for a general work on the $TV$-norm.
For an $M\times M$ array $z=(z_{s,t})_{s=1,t=1}^{M,M}$ the $TV$-norm
is defined as the convex function $TV:\mathbb{R}^{M\times M}\to\mathbb{R}$
given by 
\begin{equation}
\begin{aligned}TV(z):=\sum_{s=1}^{M-1}\sum_{t=1}^{M-1} & \sqrt{(z_{s,t}-z_{s+1,t})^{2}+(z_{s,t}-z_{s,t+1})^{2}}\\
 & \quad+\sum_{s=1}^{M-1}|z_{s,M}-z_{s+1,M}|+\sum_{t=1}^{M-1}|z_{M,t}-z_{M,t+1}|.
\end{aligned}
\end{equation}

Choosing the $TV$-norm as the objective function in the ``$x$-space''
or the ``$y$-space'', or both, and associating it with the feasibility-seeking
problem of the fully-discretized inverse problem of IMRT (Problem
\ref{prob:IMRT0}) leads naturally to formulations of the SMP in Problem
\ref{prob:smp-1}.

\subsubsection{The experimental setup}

For the purpose of our numerical experiment we confine ourselves specifically
to a case of the feasibility-seeking problem of the fully-discretized
inverse problem of IMRT (Problem \ref{prob:IMRT0}) where there are
$L$ tumor structures and the whole rest of the cross-section is considered
as one single organ-at-risk, i.e., we let $S_{r}=S\subseteq\{1,2,\ldots,m\}$
for all $r=1,2,\ldots,R.$ This leads to the next split problem of
minimizing the $TV$-norm of the dose subvectors so that uniformity
of dose distribution will be achieved for each tumor structure separately.
\begin{problem}
\label{prob:IMRT} Let $\mathbb{R}^{n}$ be the $x$-space of intensity
vectors, let $\mathbb{R}^{m}$ be the $y$-space of dose vectors,
and $A$ be the dose matrix relating them to each other. For $\ell=1,2,\ldots,L$
denote by $T_{\ell}\subseteq\{1,2,\ldots,m\}$ the sets of pixels
corresponding to the $\ell$-th tumor structure and let $S\subseteq\{1,2,\ldots,m\}$
be the complementary set of pixels that do not belong to any of the
target structures and represent all organs at risk. Set $0$ and $\overline{e}$
the lower and upper bounds for the beamlets intensities. Let $\underline{d}$
and $\overline{d}$, and $\underline{c}_{\ell}$ and $\overline{c}_{\ell}$
be the dose bounds for pixels in an organ-at-risk and at the tumor
structures, respectively. We further assume that the dose vector $y=Ax$
consists of $L+1$ subvectors $y=(y^{\ell})_{\ell=1}^{L+1}$such that
the first $L$ subvectors consists of the doses absorbed in pixels
of the $L$ tumor structures and $y^{L+1}$ is the dose absorbed in
the complementary tissue $S.$

The task is to find an intensities vector $x$ such that 
\begin{equation}
\begin{aligned}\underline{c}_{\ell}\leq y_{j}^{\ell}=\sum_{i=1}^{n}a_{ij}x_{i}\leq\overline{c}_{\ell}, & \quad\text{ for all }j\in T_{\ell},\quad\ell\in{\{1,2,\ldots,L\}},\\
\underline{d}\leq y_{j}^{L+1}=\sum_{i=1}^{n}a_{ij}x_{i}\leq\overline{d}, & \quad\text{ for all }j\in S,\\
\underline{e}\leq x_{i}\leq\overline{e}, & \quad\text{ for all }i\in{\{1,2,\ldots,n\}},
\end{aligned}
\end{equation}
\begin{equation}
\text{ and }y^{\ell}\in\argmin\{TV(u)\:|\:u\in[\underline{c}_{\ell},\overline{c}_{\ell}]^{|T_{\ell}|}\}\text{ for all }\ell\in{\{1,2,\ldots,L\}},
\end{equation}
where, for every $\ell\in{\{1,2,\ldots,L\}}$, $y^{\ell}$ denotes
the subvector of the vector $y$ associated with the $\ell$-th tumor,
and $|T_{\ell}|$ is the cardinality of the set $T_{\ell}$.
\end{problem}

\begin{figure}[H]
\centering \includegraphics[height=0.35\textwidth]{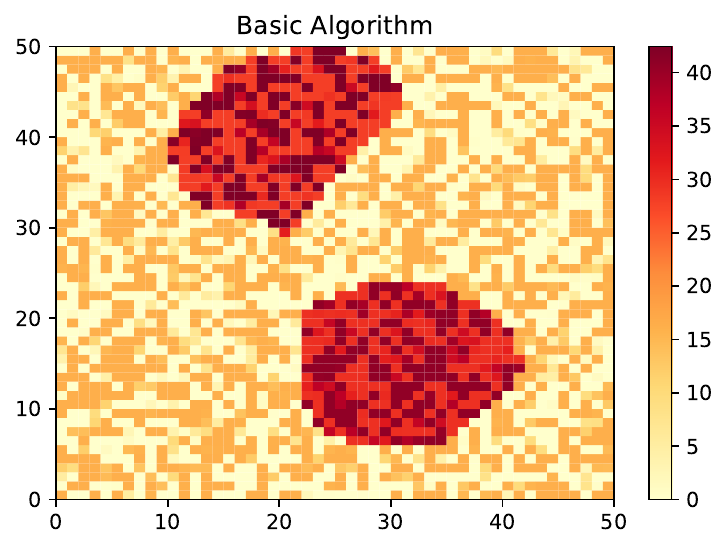}
\includegraphics[height=0.35\textwidth]{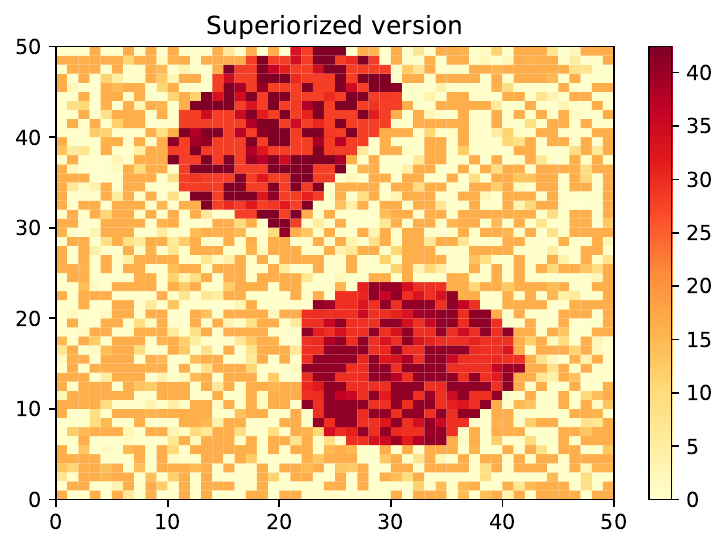}
\includegraphics[height=0.35\textwidth]{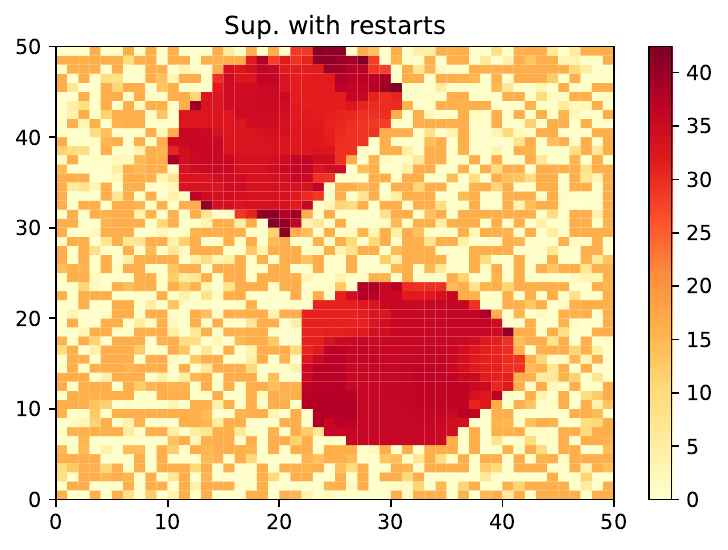}
\caption{Heat maps of the solutions in the $y$-space of pixel doses for the
first run of each one of the studied algorithms. Represented is a
$50\times50$ square grid of pixels, where the color indicates the
dose absorbed in each pixel. \label{fig:heatmap}}
\end{figure}

This is the problem we worked on in our experiment. We do not use
real data but replicate a realistic situation. In particular, we consider
a cross-section of $50\times50$ square pixels, which translates into
the dose vector in the $y$-space $\mathbb{R}^{2500}$. The number
of external radiation beamlets is $n=2840$, meaning that the $x$-space
is $\mathbb{R}^{2840}$. In the cross-section we have two tumor structures
of irregular shapes, whose location appears in Figure~\ref{fig:heatmap}.
In order to guarantee the existence of a feasible solution for Problem~\ref{prob:IMRT},
we generated the data as follows.
\begin{itemize}
\item We generate a vector $\overline{y}\in\mathbb{R}^{2500}$ with components
randomly distributed in the interval $[0,15]$ for the pixels corresponding
to organs-at-risk, and in the interval $[10,40]$ for pixels of tumor
structures.
\item We randomly generated a matrix $V\in\mathbb{R}^{2840\times2500}$
with entries in the interval $[0,1]$ and defined the dose matrix
$A\in\mathbb{R}^{2500\times2840}$, mapping the $x$-space onto the
$y$-space, as the generalized left inverse of $V$, i.e., we took
$A:=(V^{T}V)^{-1}V^{T}$.
\item We defined $\overline{x}:=V\overline{y}$, which implies that $\overline{y}=A\overline{x}$.
\item We set the bounds for the constraints of~Problem \ref{prob:IMRT}
as
\begin{equation}
\left\{ \begin{array}{l}
\underline{d}=0,\,\overline{d}=\max\{\overline{y}_{j}\mid j\in S\}+5\varepsilon_{1},\\
\underline{c}_{1}=\min\{\overline{y}_{j}\mid j\in T_{1}\}-5\varepsilon_{2},\,\overline{c}_{1}=\max\{\overline{y}_{j}\mid j\in T_{1}\}+5\varepsilon_{3},\\
\underline{c}_{2}=\min\{\overline{y}_{j}\mid j\in T_{2}\}-5\varepsilon_{4},\,\overline{c}_{2}=\max\{\overline{y}_{j}\mid j\in T_{2}\}+5\varepsilon_{5},\\
\underline{e}=(\varepsilon_{6}+1)/2\min\{\overline{x}_{i}\mid i\in{\{1,2,\ldots,n\}}\},\\
\overline{e}=(1+\varepsilon_{7}/2)\max\{\overline{x}_{i}\mid i\in{\{1,2,\ldots,n\}}\},
\end{array}\right.
\end{equation}
where the sub-indices in $\overline{c}$ and $\underline{c}$ refer
to the first and second tumor structures and, for $i\in{\{1,2,\ldots,7\}}$,
$\varepsilon_{i}$ are randomly picked real numbers in the interval
$(0,1]$.
\end{itemize}
These choices during the data generation guarantee that there exists
a feasible solution for Problem~\ref{prob:IMRT} with these data,
namely $\overline{x}$.

In our experimental work, we ran the basic algorithm (Algorithm~\ref{a:BA_svip})
and the superiorized version of the basic algorithm (Algorithm~\ref{a:SV_SPM})
with and without restarts. For all of them we took the algorithmic
operator as 
\begin{equation}
\mathbf{T}:=P_{\textbf{V}}\circ\left(P_{[\underline{e},\overline{e}]}\times P_{Q}\right),
\end{equation}
with $\mathbf{V}:=\left\{ (x,y)\in\mathbb{R}^{n}\times\mathbb{R}^{m}\mid Ax=y\right\} $
and $P_{Q}:\mathbb{R}^{m}\rightarrow\mathbb{R}^{m}$ defined component-wise
as 
\begin{equation}
P_{Q}(y_{j}):=\left\{ \begin{array}{l}
P_{[\underline{d},\overline{d}]},\text{ if }j\in S,\\
P_{[\underline{c}_{1},\overline{c}_{1}]},\text{ if }j\in T_{1},\\
P_{[\underline{c}_{2},\overline{c}_{2}]},\text{ if }j\in T_{2}.
\end{array}\right.
\end{equation}
We tested the three algorithms with different choices of the parameters
and present here the most advantageous for each one.
Specifically, in Algorithm \ref{a:SV_SPM} with or without restart
the step-sizes were taken in the sequence  $\{c\,\alpha^{\ell}\}_{\ell=0}^{\infty}$
with a constant kernel $\alpha$ and a positive number $c$, and we took $N=5$. We performed some experiments in order to determine the best choice of $\alpha$
and $c$ for each method. The results are shown in Table~\ref{table:par}.
We chose $\alpha=0.999$ and $c=100\,000$ for the superiorized algorithm,
since these parameters provide the best reduction in $TV$-norm values
while performing the fastest. For the superiorization with
restarts, all of the combinations of parameters, except of the first one,  provide a  great reduction in the $TV$-norm values with respect to superiorization with no restarts. Among these combinations, $\alpha=0.99$
and $c=100$ was the fastest, so we opted for it.

The target functions $\phi_{b}$ were always the appropriate $TV$-norms.
Since no smoothing of the intensities vectors is included in the experiment,
we took $v_{x}^{k,j}=0,$ for all $k$ and $j$. The
final parameters of the two methods are the following:
\begin{itemize}
\item \textbf{Superiorization:} We took $\alpha=0.999$, $c=100\,000$
and $N=5$, and $v_{y}^{k,b,j}$ was defined as the
nonascending direction given by Theorem~\ref{t:nonasc}.
\item \textbf{Superiorization with restarts:} We took $\alpha=0.99$, $c=100$,
$W_{r}=20$ for all $r$ and $N=5$, and $v_{y}^{k,b,j}$ was defined
as the nonascending direction given by Theorem~\ref{t:nonasc}.
\end{itemize}

\begin{table}[H]
\resizebox{\columnwidth}{!}{%
\begin{tabular}{ccrrrrr|rrrrr}
\cline{3-12} \cline{4-12} \cline{5-12} \cline{6-12} \cline{7-12} \cline{8-12} \cline{9-12} \cline{10-12} \cline{11-12} \cline{12-12}
 &  & \multicolumn{5}{c|}{$\alpha=0.99$} & \multicolumn{5}{c}{$\alpha=0.999$}\tabularnewline
\cline{3-12} \cline{4-12} \cline{5-12} \cline{6-12} \cline{7-12} \cline{8-12} \cline{9-12} \cline{10-12} \cline{11-12} \cline{12-12}
\multirow{1}{*}{} &  & $c=10$ & $c=100$ & $c=1000$ & $c=10\,000$ & $c=100\,000$ & $c=10$ & $c=100$ & $c=1000$ & $c=10\,000$ & $c=100\,000$\tabularnewline
\hline
\multirow{3}{*}{Sup.} & $TV1$ & 2433.38 & 2433.36 & 2432.93 & 2430.40 & 2430.40 & 2432.26 & 2423.20 & 2300.88 & 2167.01 & 2166.97\tabularnewline
 & $TV2$ & 3056.41 & 3056.34 & 3055.51 & 3052.28 & 3052.28 & 3054.67 & 3039.15 & 2899.47 & 2714.62 & 2714.58\tabularnewline
 & Time & 204.34 & 202.87 & 203.46 & 201.12 & 208.12 & 202.69 & 201.84 & 199.67 & 199.33 & 196.57\tabularnewline
\hline
\multirow{3}{*}{Sup. Restarts} & $TV1$ & 2078.81 & 368.81 & 397.35 & 381.65 & 371.03 & 249.96 & 361.70 & 386.13 & 381.51 & 382.51\tabularnewline
 & $TV2$ & 2689.98 & 707.2 & 598.24 & 579.47 & 568.02 & 585.33 & 583.39 & 534.21 & 528.54 & 528.56\tabularnewline
 & Time & 367.65 & 453.47 & 598.55 & 779.16 & 950.63 & 2362.15 & 3560.19 & 4426.75 & 6181.23 & 7174.09\tabularnewline
\hline
\end{tabular}}\caption{Average $TV$-norm values for the first and second subvectors and
average time (in seconds) obtained by running the Superiorized and
Superiorized with restarts algorithm with different choices of parameters
for 10 random initial points. The algorithms were stopped
when a proximity of 0.01 was reached.}\label{table:par}
\end{table}
We performed multiple runs of the three algorithms. At each run, each
of the algorithms was initialized at the same starting point which
was randomly generated in the interval $[\underline{e},\overline{e}]$.
We define the proximity of an iterate as the distance to the feasible
region, i.e., for an iterate pair $(x^{k},y^{k})$, we define its
proximity as 
\begin{equation}
{\rm prox}(x^{k},y^{k}):=\|x^{k}-P_{[\underline{e},\bar{e}]}(x^{k})\|+\|y^{k}-P_{Q}(y^{k})\|.
\end{equation}

Note that, due to the definition of the algorithmic operator $\mathbf{T}$,
the distance of $(x^{k},y^{k})$ to $\mathbf{V}$ is $0$. All three
algorithms were terminated once the proximity became less than $0.01$.
The obtained results for all different runs are summarized in Table~\ref{tab:1}.
Our numerical experiments showed that superiorization with restarts
was considerably the best performing algorithm regarding the target
function reduction, while superiorization alone, without restarts,
did not achieve a significant reduction with respect to the basic
algorithm.

This fact can be graphically observed in the heat maps of Figure~\ref{fig:heatmap},
where we represent the dose in the pixels of the cross-section at
the last iteration of each algorithm. The uniformity of the heat distribution
in a tumor structure represents the dose distribution in that structure.
Clearly, superiorization with restarts provided a more homogeneous
dose distribution in the tumorous pixels. We observed the increased
uniformity of dose distributions in the tumors in all our algorithmic
runs of the superiorization with restarts method. However, depending
on the datasets and the allowable parameters the level of the uniformity
may vary.

\begin{figure}[H]
\centering \includegraphics[width=0.49\textwidth]{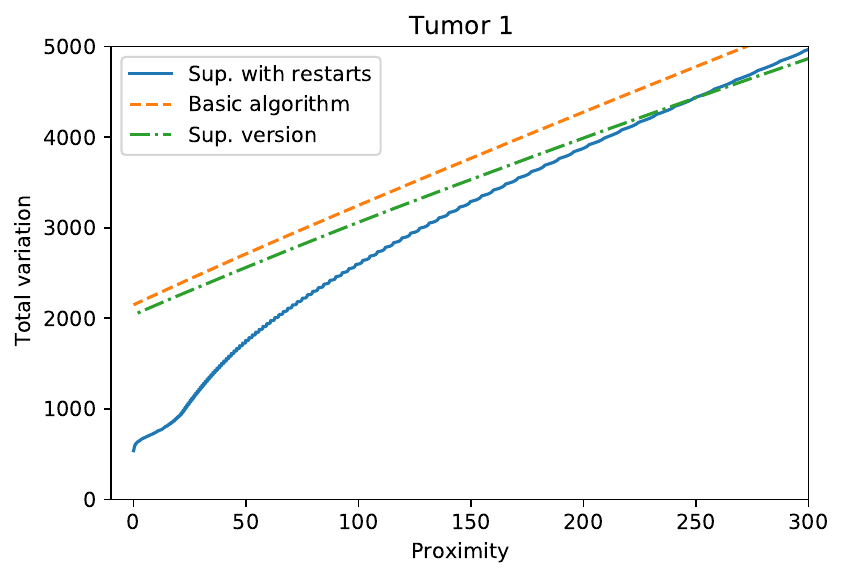}
\includegraphics[width=0.49\textwidth]{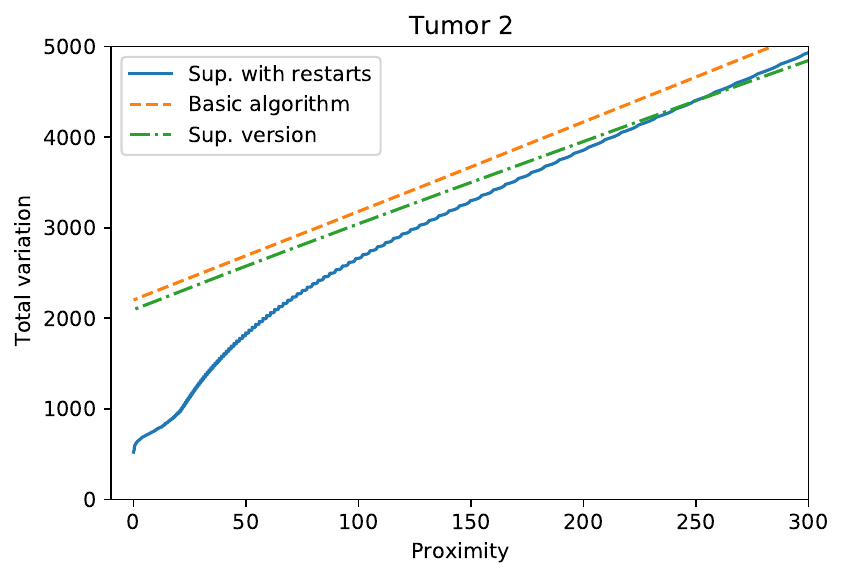}
\caption{The evolution of the total variation and the proximity of the iterations
of the first run of each of the algorithms for the subvector associated
to the first tumor (left) and the second tumor (right). In these \textquotedblleft proximity-target
function curves\textquotedblright{} the iteration indices $k$ increase
from right to left in each of the plots. \label{fig:proximityvstv}}
\end{figure}

The evolution along the iterations of the proximity and the total
variation of the algorithms is shown in Figure~\ref{fig:proximityvstv}
with ``proximity-target function curves'' (which were introduced
in \cite{DFS-2021}), where the iteration indices $k$ increase from
right to left in each of the plots. Finally, we note that superiorization
with restarts needed more time and a larger number of iterations to
reach the desired proximity.

In our experiments we have observed that other choices of parameters
for the superiorization with restarts runs can be employed to reduce
its running times and make them comparable to those of the superiorized
algorithm without restarts and, at the same time, still achieve a
significant reduction of the target function when compared to the
other algorithms.

\begin{table}[H]
\begin{centering}
\begin{tabular}{ccr@{\extracolsep{0pt}.}lr@{\extracolsep{0pt}.}lr@{\extracolsep{0pt}.}lr@{\extracolsep{0pt}.}lr@{\extracolsep{0pt}.}l}
\toprule
Run  &  & \multicolumn{2}{c}{1} & \multicolumn{2}{c}{2} & \multicolumn{2}{c}{3} & \multicolumn{2}{c}{4} & \multicolumn{2}{c}{5}\tabularnewline
\midrule
\multirow{3}{*}{$TV$ for subvector 1} & Basic  & 2405&3  & 2498&4  & 2289&3  & 2624&74 & 2474&9\tabularnewline
 & Superiorized  & 2072&8 & 2230&1  & 2089&6  & 2362&7  & 227&2\tabularnewline
 & Sup. Restarts  & 421&9  & 315&6  & 404&3  & 349&3  & 301&1\tabularnewline
\midrule
\multirow{3}{*}{$TV$ for subvector 2} & Basic  &  3019&1 & 3252&8  & 3002&5  & 3076&6  & 3096&3\tabularnewline
 & Superiorized  & 2703&3  & 2837&4  & 2558&9  & 2744&3  & 2848&2\tabularnewline
 & Sup. Restarts  & 817&7  & 759&7  & 688&1  & 709&2  & 531&2\tabularnewline
\midrule
\multirow{3}{*}{Run time (sec.)} & Basic  & 96&1  & 94&6  & 91&6 & 93&6 & 95&4\tabularnewline
 & Superiorized  & 261&6  & 263&3  & 253&8  & 269&4  & 267&3\tabularnewline
 & Sup. Restarts  & 577&5  & 588&0  & 630&9  & 580&3  & 582&7\tabularnewline
\midrule
\multirow{3}{*}{No. of iterations} & Basic  & \multicolumn{2}{c}{7352} & \multicolumn{2}{c}{7265} & \multicolumn{2}{c}{7117} & \multicolumn{2}{c}{7400} & \multicolumn{2}{c}{7505}\tabularnewline
 & Superiorized  & \multicolumn{2}{c}{7327} & \multicolumn{2}{c}{7195} & \multicolumn{2}{c}{7095} & \multicolumn{2}{c}{7382} & \multicolumn{2}{c}{7478}\tabularnewline
 & Sup. Restarts  & \multicolumn{2}{c}{14880} & \multicolumn{2}{c}{14819} & \multicolumn{2}{c}{17140} & \multicolumn{2}{c}{15275} & \multicolumn{2}{c}{14778}\tabularnewline
\cmidrule{1-2}\cmidrule{2-12} \cmidrule{3-12} \cmidrule{5-12} \cmidrule{7-12} \cmidrule{9-12} \cmidrule{11-12}
\end{tabular}
\par\end{centering}
\centering{}\caption{$TV$-norm values for the first and second subvector, run times and
number of iterations resulting from running the Basic, Superiorized
and Superiorized with restarts algorithms (for runs with 5 different
random initial points) until a proximity of $0.01$ was reached\label{tab:1}.}
\end{table}

\section*{Acknowledgments}

We thank the referee for his constructive comments and criticisms
that helped to improve our work. Yair Censor gratefully acknowledges
enlightening discussions with Professor Reinhard Schulte from Loma
Linda University in Loma Linda, California, about minimization or
superiorization on subvectors in IMRT treatment planning. The work
of Yair Censor was supported by the ISF-NSFC joint research plan Grant
Number 2874/19. Francisco Arag\'on and David Torregrosa were partially
supported by the Ministry of Science, Innovation and Universities
of Spain and the European Regional Development Fund (ERDF) of the
European Commission, Grant PGC2018-097960-B-C22, and the Generalitat
Valenciana (AICO/2021/165). David Torregrosa was supported by MINECO
and European Social Fund (PRE2019-090751) under the program ``Ayudas
para contratos predoctorales para la formaci\'on de doctores'' 2019.

\end{document}